\documentclass[11pt,a4paper]{article}
\usepackage{epsfig}
\usepackage{graphics}
\usepackage{amsmath,amsfonts,bbm,amssymb}
\usepackage{subfigure}
\usepackage{graphicx}
\usepackage{epstopdf}
\usepackage{hyperref}
\usepackage{pdfsync}
\usepackage{color}
\usepackage{tikz}
\usetikzlibrary{patterns}

\newtheorem{thm}{Theorem}

\newtheorem{proposition}{Proposition}
\newtheorem{definition}{Definition}
\newtheorem{remark}{Remark}

\numberwithin{equation}{section}
\numberwithin{thm}{section}
\numberwithin{lemma}{section}
\numberwithin{proposition}{section}
\numberwithin{definition}{section}
\numberwithin{remark}{section}

\def\arginf{\mathop{\hbox{\rm arginf }}}
\def\Alpha{\boldsymbol{\alpha}}
\def\Beta{\boldsymbol{\beta}}
\def\Zeta{\boldsymbol{\zeta}}

\title{A unified treatment for non-asymptotic and asymptotic approaches to minimax signal detection\footnote{This work is dedicated to the memory of our collaborator and friend Yuri I. Ingster who passed away on August 8th, 2012.}}
\author{{\em Cl\'ement ~Marteau}, \\
         Institut Math\'ematiques de Toulouse,
         INSA de Toulouse,\\
         Universit\'e de Toulouse,\\
         135, avenue de Rangueil, 31 077 Toulouse Cedex 4, France.\\
        {\em Email}:~\texttt{clement.marteau@math.univ-toulouse.fr}
          \\ \\
        and
          \\ \\
        {\em Theofanis ~Sapatinas},\\
        Department of Mathematics and Statistics,\\
        University of Cyprus,\\
        P.O. Box 20537,
        CY 1678 Nicosia,
        Cyprus.\\
        {\em Email}:~\texttt{fanis@ucy.ac.cy}}

\begin{document}
\maketitle

\begin{abstract}
We are concerned with minimax signal detection. In this setting, we discuss  non-asymptotic and asymptotic approaches through a unified treatment. In particular, we consider a Gaussian sequence model that contains classical models as special cases, such as, direct, well-posed inverse and ill-posed inverse problems. Working with certain ellipsoids in the space of squared-summable sequences of real numbers, with a ball of positive radius removed, we compare the construction of lower and upper bounds for the minimax separation radius (non-asymptotic approach) and the minimax separation rate (asymptotic approach) that have been proposed in the literature. Some additional contributions, bringing into light links between non-asymptotic and asymptotic approaches to minimax signal, are also presented. An example of a mildly ill-posed inverse problem is used for illustrative purposes. In particular, it is shown that tools used to derive `asymptotic' results can be exploited to draw `non-asymptotic' conclusions, and vice-versa.

% In order to enhance our understanding of these two minimax signal detection paradigms, we bring into light hitherto unknown similarities and links between non-asymptotic and asymptotic approaches. 

\medskip
\noindent
{\bf AMS 2000 subject classifications:} 62G05, 62K20\\

\medskip
\noindent
{\bf Keywords and phrases:} Gaussian sequence models, ill-posed and well-posed inverse problems, minimax signal detection.
\end{abstract}

\newpage

\section{Introduction}
We consider the following Gaussian sequence model (GSM),
\begin{equation}\label{1.0.0}
Y_j=b_j \theta_j+\varepsilon \, \xi_j, \quad j \in \mathbb{N},
\end{equation}
where $\mathbb{N} =\{1,2,\ldots \}$ is the set of natural numbers, $b=\{b_j\}_{j \in \mathbb{N}} > 0$ is a known sequence,  $\theta=\{\theta_j\}_{j\in\mathbb{N}} \in l^2(\mathbb{N})$ is the unknown signal of interest, $\xi=\{\xi_j\}_{j \in\mathbb{N}}$ is a sequence of independent standard Gaussian random variables, and $\varepsilon >0$ is a known parameter (the noise level). The observations are given by the sequence $Y=\{Y_j\}_{j \in \mathbb{N}}$ from the GSM (\ref{1.0.0}) and their joint law is denoted by $\mathbb{P}_{\theta}$. Here, $l^2(\mathbb{N})$ denotes the space of squared-summable sequence of real numbers, i.e., 
$$
l^2(\mathbb{N}) =\left\lbrace \theta \in \mathbb{R}^\mathbb{N}:\; \|\theta\|^2:=\sum_{j \in \mathbb{N}}\theta_j^2 < +\infty \right\rbrace. 
$$
The GSM (\ref{1.0.0}) arises in many well-known situations. For instance, consider  the Gaussian white noise model (GWNM)
\begin{equation}\label{1.0.2}
dX_{\varepsilon}(t)= Af(t)dt + \varepsilon \, dW(t), \quad t \in V,
\end{equation}
where $A$ is a known linear operator acting on a Hilbert
space ${\cal H}_1$  with values on another Hilbert space ${\cal H}_2$, $f(\cdot) \in {\cal H}_1$ is the unknown response function that one wants to detect or estimate, $W(\cdot)$ is a standard Wiener process on
$V \subseteq \mathbb{R}$ and $\varepsilon>0$ is a known parameter (the noise level). For the sake of simplicity, we only consider the case when $A$
is injective (meaning that $A$ has a trivial nullspace) and assume that $V=[0, 1]$, ${\cal H}_1 = L^2(V)$, $U \subseteq \mathbb{R}$ and ${\cal H}_2=L^2(U)$.

\begin{itemize}
\item (direct problem) Let $A=I$ (the identity operator). Let $\{\phi_j\}_{j \in \mathbb{N}}$ be an orthonormal basis on $L^2(V)$. Transforming the GWNM (\ref{1.0.2}) with $A=I$ into the Fourier domain, the GSM (\ref{1.0.0}) arises with $Y_j = \int_{0}^{1} \phi_j(t) dX_{\varepsilon}(t)$, $\theta_j =\int_{0}^{1} \phi_j(t)f(t)dt$, $\xi_j=\int_{0}^{1} \phi_j(t)dW(t)$ and $b_j=1$, for all $j \in \mathbb{N}$.   

\item (well-posed inverse problem) Let $A$ be a self-adjoint operator that admits an eigenvalue-eigenfunction decomposition $(b_j,\varphi_j)_{j\in \mathbb{N}}$, in the sense that
$$
A \varphi_j=b_j\varphi_j,  \quad j \in \mathbb{N},
$$
where $b_j > b_0$, for some $b_0 >0$, for all $j \in \mathbb{N}$. Thus, the GSM (\ref{1.0.0}) arises with $Y_j = \int_{0}^{1} \varphi_j(t) dX_{\varepsilon}(t)$, $\theta_j =\int_{0}^{1} \varphi_j(t)f(t)dt$, $\xi_j=\int_{0}^{1} \varphi_j(t)dW(t)$ and $b_j>b_0>0$, for all $j \in \mathbb{N}$. In this case, the GWNM (\ref{1.0.2}) corresponds to a so-called well-posed inverse problem. Possible examples of such decompositions arise
with, e.g., differential or Sturm-Liouville operators.
\item (ill-posed inverse problem) In most cases of interest, however, $A$ is a compact operator (see, e.g., Chapter 2 of \cite{EHN_1996}). In particular, it admits a singular value decomposition (SVD) $(b_j, \psi_j, \varphi_j)_{j\in \mathbb{N}}$, in the sense that
$$
A \varphi_j=b_j\psi_j, \quad A^{\star}\psi_j=b_j \varphi_j, \quad j \in \mathbb{N},
$$
where $A^{\star}$ denotes the adjoint operator of $A$ -- note that $(b_j^2)_{j \in \mathbb{N}}$ and $(\varphi_j)_{j \in \mathbb{N}}$ are, respectively,  the eigenvalues and the eigenfunctions of $A^\star A$. Thus, the GSM (\ref{1.0.0}) arises with $Y_j = \int_{0}^{1} \psi_j(t) dX_{\varepsilon}(t)$, $\theta_j =\int_{0}^{1} \varphi_j(t)f(t)dt$, $\xi_j=\int_{0}^{1} \psi_j(t)dW(t)$ and $b_j>0$ (since $A$ is injective), for all $j \in \mathbb{N}$. In this case, the GWNM (\ref{1.0.2}) corresponds to a so-called ill-posed inverse problem since the inversion of $A^*A$ is not bounded. Possible examples of such decompositions arise
with, e.g., convolution or Radon-transform operators. The effect of the ill-posedeness of the model is clearly seen in the decay of (the singular values) $b_j$ towards 0 as $j \to +\infty$. As $j \to +\infty$, $b_j \theta_j$ gets weaker and it is then more difficult to estimate or detect the sequence $\theta=\{ \theta_j\}_{j\in\mathbb{N}}$.
\end{itemize}

From the above discussion, it is evident that one can undertake statistical inference based on observations either from the GSM (\ref{1.0.0}) or from the (equivalent) GWNM (\ref{1.0.2}). Estimation in either models has received much attention over the last decades, providing also optimality results (in the minimax sense) over various loss functions and sequence/function spaces. Many methods have been considered including kernel, local polynomial, spline, projection and wavelet methods (see, e.g., \cite{Wahba_1990}, \cite{FG_1996}, \cite{Tsybakov_2009}, \cite{Cavalier_book}, \cite{Johnstone_2013}). \\

On the other hand, signal detection has received less attention. Minimax signal detection in the GSM (\ref{1.0.0}) with $b_j=1$ for all $j \in \mathbb{N}$ has been studied in \cite{E_1991} and in detail in the seminal work of \cite{IngsterI}, \cite{IngsterII} and \cite{IngsterIII} (see also \cite{IS_2003}). This work uses an asymptotic framework, that is, the noise level $\varepsilon >0$ is allowed to {\em converge} to zero. A corresponding non-asymptotic framework, that is, for any {\em fixed} value of the noise level $\varepsilon >0$, has been studied in \cite{Baraud} and \cite{Laurent_2002}. Non-asymptotic and asymptotic studies for minimax signal detection in the GSM (\ref{1.0.0}) with $b_j >0$ for all $j \in \mathbb{N}$ have been recently considered in \cite{ISS_2012} and \cite{LLM_2012}, respectively, in order to study minimax signal detection in ill-posed inverse problems. Despite the fact that the considered minimax signal detection problem is the same in the aforementioned studies, the final aims and the methodologies involved sometimes differ. \\

Bearing in mind the different issues and tasks involved, our aim below is to provide a unified treatment for non-asymptotic and asymptotic approaches to minimax signal detection in the GSM (\ref{1.0.0}). In particular, we look for common ground between them that will enhance our understanding of these two existing minimax signal detection paradigms. This paper is organized as follows.  Section \ref{msd} considers minimax signal detection from both non-asymptotic and asymptotic point of views.  Section \ref{s:control} discusses the construction of upper and lower bounds of the minimax separation radius (non-asymptotic approach) and the minimax separation rate (asymptotic approach) in a unified treatment, and points out several similarities. Section \ref{connection} brings into light hitherto unknown links between non-asymptotic and asymptotic approaches to minimax signal detection. It also contains an example of a mildly ill-posed inverse problem for illustrative purposes.  %It demonstrates that non-asymptotic and asymptotic approaches to minimax signal detection, somehow, merge. 
In particular, it is shown that tools used to derive `asymptotic' results can be exploited to draw `non-asymptotic' conclusions, and vice-versa. Finally, Section \ref{s:discussion} draws some concluding remarks and provides an avenue for future research.\\

Throughout the paper, we set the following notations.  For all $x,y\in \mathbb{R}$, $\delta_x(y)=1$ if $x=y$ and $\delta_x(y)=0$ if $x\not =y$. Also, $x \wedge y:= \inf\{x,y\}$ and $(x)_+:=\max\{0,x\}$. Given two collections $(c_\varepsilon)_{\varepsilon >0}$ and $(d_\varepsilon)_{\varepsilon >0}$ of real numbers, $c_\varepsilon \sim d_\varepsilon$ means that there exist $0<\kappa_0 \leq \kappa_1 <\infty$ such that $\kappa_0 \leq c_\varepsilon/d_\varepsilon \leq \kappa_1$ for all $\varepsilon >0$.
In the same spirit, given two sequences $(c_j)_{j \in \mathbb{N}}$ and $(d_j)_{j \in \mathbb{N}}$ of real numbers, $c_j \asymp d_j$ means that there exist $0<\kappa_0 \leq \kappa_1 <\infty$ such that $\kappa_0 \leq c_j/d_j \leq \kappa_1$ for all $j \in \mathbb{N}$. Finally, the abbreviation $o_\varepsilon(1)$ (resp. $\mathcal{O}_\varepsilon(1)$) will refer to a collection tending to 0 (resp. bounded) as $\varepsilon$ tends to $0$.  When the dependence is not explicitly required on the noise level $\varepsilon >0$, it will be simply denoted by $o(1)$ (resp. $\mathcal{O}(1)$).

\section{Minimax signal detection}
\label{msd}

Statistical estimation is concerned with a quantitative question.  Instead, we address below a qualitative question: given observations from the GSM (\ref{1.0.0}), our aim is to compare the underlying (unknown) signal $\theta \in l^2(\mathbb{N})$ to a (known) benchmark signal $\theta_0$, i.e., to test 
\begin{equation}
H_0: \theta=\theta_0 \;\; \mathrm{versus} \;\; \ H_1: \theta-\theta_0 \in \mathcal{F},
\label{testing_pb0}
\end{equation}
for some given $\theta_0$ and a given subspace $\mathcal{F}$. The statistical setting (\ref{testing_pb0}) is known as goodness-of-fit testing when $\theta_0 \neq 0$ and as signal detection when  $\theta_0 =0$.

\begin{remark}
\label{equiv_test}
{\rm Given observations from the GWNM (\ref{1.0.2}), the test (\ref{testing_pb0}) is related to the test
\begin{equation}
H_0: f=f_0 \;\; \mathrm{versus} \;\; H_1: f-f_0 \in \tilde{\mathcal{F}},
\label{testing_pbf}
\end{equation}
for a given benchmark function $f_0$ and a given subspace $\tilde{\mathcal{F}}$. In most cases, $\tilde{\mathcal{F}}$ contains functions $f \in L^2(V)$ that admit a Fourier series expansion with Fourier coefficients $\theta$ belonging to ${\cal F}$ (see, e.g., \cite{IS_2003}, Section 3.2). In these cases, the problems (\ref{testing_pb0}) and (\ref{testing_pbf}) are equivalent.} 
\end{remark}

The choice of the set $\mathcal{F}$ is important. Indeed, it should be rich enough in order to contain the true $\theta$. At the same time, if it is too rich, it will not be possible to control the performances of a given test due to the complexity of the problem. The common approach for such problems is to impose both a {\em regularity} condition (which characterizes the smoothness of the underlying signal) and an {\em energy} condition (which measures the amount of the underlying signal).  \\ 

Concerning the regularity condition, we will work with certain ellipsoids in $l^2(\mathbb{N})$. In particular, we assume that $\theta\in \mathcal{E}_{a}(R)$, the set $\mathcal{E}_{a}(R)$ being defined as
$$ 
\mathcal{E}_{a}(R) = \left\lbrace \theta\in l^2(\mathbb{N}), \ \sum_{j \in \mathbb{N}} a_j^2 \theta_j^2 \leq R \right\rbrace,
$$
where $a=(a_j)_{j\in \mathbb{N}}$ denotes a non-decreasing sequence of positive real numbers with $a_j \rightarrow +\infty$ as $j \rightarrow +\infty$,
and $R>0$ is a constant. The set $ \mathcal{E}_{a}(R)$ can be seen as a condition on the decay of $\theta$. The cases where $a$ increases very fast correspond to $\theta$ with a small amount of non-zero coefficients. In such a case, the corresponding signal can be considered as being `smooth'. \\

Without loss of generality, in what follows, we set $R=1$. In order to simplify the notation, we will avoid the dependency to this term in all related quantities. In particular, we will write $\mathcal{E}_{a}$ instead of $\mathcal{E}_{a}(1)$. \\

Regarding the energy condition, it will be measured in the $l^2(\mathbb{N})$-norm. In particular, given $r_\varepsilon>0$ (called the radius), which is allowed to depend on the noise level $\varepsilon >0$, we will consider $\theta\in \mathcal{E}_{a}$ such that $\| \theta\| > r_\varepsilon$.  Given a smoothness sequence $a$ and a radius $r_\varepsilon>0$, the set $\mathcal{F}$ can thus be defined as
\begin{equation}
\label{def:f-set}
\mathcal{F} := \Theta_a(r_\varepsilon) = \left\lbrace \theta \in \mathcal{E}_a, \ \|\theta \| \geq r_\varepsilon \right\rbrace.
\end{equation}
%This means that the set ${\cal F}$ is an ellipsoid in $l^2(\mathbb{N})$ with  a ball of radius $r_\varepsilon >0$ removed. We stress that several investigations have been conducted in the case where the energy \textcolor{red}{$\|\theta\|$} is measured with a different norm (e.g., $l^p(\mathbb{N})$ with $p\not = 2$). We can mention, for instance, \cite{Baraud} or \cite{IS_2003}, among others. However, for the sake of clarity, this setting will not be considered in this survey.   \\
Since $\theta_0$ and $b_j>0$, $j\in\mathbb{N}$, are known, and assuming that $\theta_0 \in \mathcal{E}_a$, without loss of generality, given observations from the GSM (\ref{1.0.0}), we restrict ourselves to the hypothesis testing setting (\ref{testing_pb0}) with $\theta_0=0$ (i.e., signal detection).\\

In summary, given observations from the GSM (\ref{1.0.0}), we will be dealing with the following signal detection problem
\begin{equation}
H_0: \theta=0 \ \mathrm{versus} \ H_1: \theta\in \Theta_a(r_\varepsilon),
\label{testing_pb1}
\end{equation}
where $\Theta_a(r_\varepsilon)$ is defined in (\ref{def:f-set}). The sequence $a$ being fixed, the main issue for the problem (\ref{testing_pb1}) is then to characterize the values of $r_\varepsilon >0$ for which both hypotheses $H_0$ (called the null hypothesis) and $H_1$ (called the alternative hypothesis) are `separable'  (in a sense which will be made precise later on).\\

In the following, a (non-randomized) test $\Psi:=\Psi(Y)$ will be defined as a measurable function of the observation $Y=(Y_j)_{j\in\mathbb{N}}$ from the GSM (\ref{1.0.0}) having values in the set $\lbrace 0,1 \rbrace$. By convention, $H_0$ is rejected if $\Psi=1$ and  $H_0$ is not rejected if $\Psi=0$. Then, given a test $\Psi$, we can investigate 
\begin{itemize}
\item the first kind error probability defined as 
\begin{equation}\Alpha_\varepsilon(\Psi):= \mathbb{P}_{0}( \Psi =1),
\label{eq:type1-F}
\end{equation} 
which measures the probability to reject $H_0$ when $H_0$ is true (i.e., $\theta=0$); it is often constrained as being bounded by a prescribed level $\alpha \in ]0,1[$, and
\item the maximal second kind error probability defined as 
\begin{equation}
\Beta_\varepsilon(\Theta_a(r_\varepsilon),\Psi) := \sup_{\theta\in \Theta_a(r_\varepsilon)} \mathbb{P}_\theta(\Psi=0),
\label{eq:type2}
\end{equation}
which measures the worst possible probability not to reject $H_0$ when $H_0$ is not true (i.e., when $\theta \in \Theta_a(r_\varepsilon))$; one would like to ensure that it is (asymptotically) bounded by a prescribed level $\beta \in ]0,1[$.
\end{itemize}
\vspace{0.2cm}

For simplicity in our exposition, we will restrict ourselves to $\alpha$-level tests, the value of $\alpha \in \,]0,1[$ being fixed.

\begin{definition}
\label{a-level-test}
A test $\Psi_\alpha$ is called an $\alpha$-level test if 
$$\Alpha_\varepsilon(\Psi_{\alpha}) \leq \alpha. $$
\end{definition}

Given the trivial test $\Psi_{\alpha}:=\alpha \in \, ]0,1[$, which does not depend on any observation, and extending the definition of a (non-randomized) test to a randomized test\footnote{a measurable function 
$\Psi:=\Psi(Y)$ of the observation $Y=(Y_j)_{j\in\mathbb{N}}$ from the GSM (\ref{1.0.0})
with values in the interval $[0,1]$: the null hypothesis is rejected with probability $\Psi(Y)$ and it is not rejected with probability $1-\Psi(Y)$. In this case, $\Alpha_\varepsilon(\Psi):= \mathbb{E}_{0}( \Psi(Y))$ and  $\Beta_\varepsilon(\Theta_a(r_\varepsilon),\Psi) := \sup_{\theta\in \Theta_a(r_\varepsilon)} \mathbb{E}_\theta(1-\Psi(Y)))$.}, it is easily seen that 
$$
\inf_{\tilde \Psi_\alpha: \,\Alpha_\varepsilon(\tilde \Psi_\alpha)\leq \alpha} \; \Beta_ {\varepsilon}(\Theta_a(r_\varepsilon), \tilde \Psi_\alpha) \in  [0,1-\alpha], \quad \text{for all}\quad \alpha \in ]0,1[ 
$$ 
(see, e.g., \cite{IS_2003}, pp. 10-11).\\

\begin{definition}
\label{trivial}
A minimax hypothesis testing problem 
$$
H_0: \theta = 0\quad \text{versus} \quad H_1: \theta \in {\cal G},
$$
for some set ${\cal G}$ (with $0 \not \in {\cal G}$),  is called trivial if 
$$ 
\inf_{\tilde \Psi_\alpha: \,\Alpha_\varepsilon(\tilde \Psi_\alpha)\leq \alpha} \; \Beta_ {\varepsilon}({\cal F}, \tilde \Psi_\alpha)=1-\alpha \quad \text{for all} \quad \alpha \in ]0,1[,
$$
and is called asymptotical trivial if 
$$ 
\inf_{\tilde \Psi_\alpha: \,\Alpha_\varepsilon(\tilde \Psi)\leq \alpha} \; \Beta_ {\varepsilon}({\cal F}, \tilde \Psi_\alpha) =1-\alpha +o_\varepsilon(1)  \quad \text{for all} \quad \alpha \in ]0,1[. 
$$ 
\end{definition}

The regularity and energy conditions imposed above, when taken together, i.e., when ${\cal F}$ is given by (\ref{def:f-set}), result (provided the radius $r_\varepsilon >0$ is `well-chosen') in a non-trivial or an  asymptotical non-trivial minimax signal detection problem (\ref{testing_pb1}).  This means, in particular, that both hypotheses $H_0$ and $H_1$ are, in some sense, separable in such a framework. Two different point of views, the so-called non-asymptotic and asymptotic minimax signal detection approaches, are at hand, that have been respectively developed in, e.g., \cite{Baraud}, \cite{Laurent_2002}, \cite{LLM_2012} and \cite{ISS_2012}. We elaborate on both approaches in the subsequent sections. \\

\subsection{The non-asymptotic approach}
Let $\alpha,\beta\in ]0,1[$ be given, and let $\Psi_\alpha$ be an $\alpha$-level test. 

\begin{definition}
The separation radius of the $\alpha$-level test $\Psi_\alpha$ over the class $\mathcal{E}_a$ is defined as
$$ r_\varepsilon(\mathcal{E}_a,\Psi_\alpha,\beta) := \inf \left\lbrace r_\varepsilon>0: \ \Beta_\varepsilon(\Theta_a(r_\varepsilon),\Psi_\alpha)  \leq \beta\right\rbrace,$$
where the maximal second kind error probability $\Beta_\varepsilon(\Theta_a(r_\varepsilon),\Psi_\alpha)$ is defined in (\ref{eq:type2}).
\label{eq:minimax_separation_radius}
\end{definition}

In some sense, the separation radius $r_\varepsilon(\mathcal{E}_a,\Psi_\alpha,\beta)$ corresponds to the smallest possible value of the available signal $\| \theta \|$ for which $H_0$ and $H_1$ can be `separated' by the $\alpha$-level test $\Psi_\alpha$ with prescribed first and maximal second kind error probabilities, $\alpha$ and $\beta$, respectively.

\begin{definition}
\label{def:minimax_radius}
The minimax separation radius $\tilde{r}_\varepsilon:=\tilde{r}_{\varepsilon}(\mathcal{E}_a, \alpha, \beta)>0$ over the class $\mathcal{E}_a$ is defined as
\begin{equation}
\tilde r_{\varepsilon}:= \inf_{\tilde \Psi_\alpha:\Alpha_\varepsilon(\tilde \Psi_\alpha)\leq \alpha} r_\varepsilon(\mathcal{E}_a, \tilde\Psi_\alpha,\beta).
\label{eq:minimax_radius}
\end{equation}
\end{definition}

The minimax separation radius $\tilde r_{\varepsilon}$ corresponds to the smallest radius $r_{\varepsilon} >0$ such that there exists some $\alpha$-level test $\tilde \Psi_\alpha$ for which the maximal second kind error probability $\Beta_\varepsilon(\Theta_a(r_\varepsilon),\tilde \Psi_\alpha)$ is not greater than $\beta$.\\
%, uniformly over all $\alpha$-level tests $\tilde \Psi_\alpha$.  

It is worth mentioning that Definitions \ref{eq:minimax_separation_radius} and \ref{def:minimax_radius} are valid for any {\em fixed} $\varepsilon >0$ (i.e., it is not required that $\varepsilon \rightarrow 0$). The performances of any given test $\Psi_\alpha$ is easy to handle in the sense that the first kind error probability $\Alpha_\varepsilon(\Psi_\alpha)$ is bounded by $\alpha$ (i.e., $\Psi_\alpha$ is an $\alpha$-level test), and that the dependence of the minimax separation radius $\tilde r_{\varepsilon}$ with respect to given $\alpha$ and $\beta$ can be precisely described. \\

In practice, given an $\alpha$-level test $\Psi_\alpha$, it might be appropriate to compare its separation radius 
$r_\varepsilon(\mathcal{E}_a,\Psi_\alpha,\beta)$ to the minimax separation radius $\tilde r_{\varepsilon}$. Hence, the following definition is in order (see, e.g., \cite{Baraud}, \cite{LLM_2012}).

\begin{definition}
\label{def:powerFanis}
A $\alpha$-level test $\Psi_\alpha$ is said to be \textit{powerful} over the class $\mathcal{E}_a$ if there exists a constant $\mathcal{C} \geq1$ such that, for all $\varepsilon>0$,
$$ \Beta_\varepsilon(\Theta_a(\mathcal{C}\tilde r_\varepsilon),\Psi_\alpha)  \leq \beta,$$
or, equivalently,
$$\ r_\varepsilon(\mathcal{E}_a,\Psi_\alpha,\beta) \leq \mathcal{C} \tilde r_{\varepsilon},$$
for any given $\beta \in ]0,1[$.
\end{definition}

According to Definition \ref{def:powerFanis}, for every $\varepsilon >0$, the separation radius $\ r_\varepsilon(\mathcal{E}_a,\Psi_\alpha,\beta)$ of a powerful $\alpha$-level test $\Psi_{\alpha}$ is of the order (up to a constant) of the minimax 
separation radius $\tilde r_{\varepsilon}$. In some sense, a powerful test appears to be {\em rate-optimal}.\\

We present below a general strategy for obtaining the minimax separation radius $\tilde r_{\varepsilon}$ (that implicitly also produces a powerful $\alpha$-level test $\Psi_{\alpha}$). Given an ellipsoid $\mathcal{E}_a$, one has to find a radius $r_\varepsilon^\star >0$ such that
$$\mbox{(Lower bound)} \quad  \tilde r_{\varepsilon} \geq r_\varepsilon^\star,$$
and to construct a specific $\alpha$-level test $\Psi_\alpha$ for which
$$ \mbox{(Upper bound)} \quad r_\varepsilon(\mathcal{E}_a,\Psi_\alpha,\beta) \leq \mathcal{C} r_\varepsilon^\star,$$ 
for some (explicitly obtained) constant $\mathcal{C} \geq 1$. It can be then easily seen that
$$ r_\varepsilon^\star \leq \tilde r_{\varepsilon} \leq \mathcal{C} r_\varepsilon^\star.$$ 

\noindent
More precisely,\\

\noindent
{\bf Lower bound:} It is enough to bound from below the following quantity
$$
\inf_{\tilde \Psi_\alpha: \,\Alpha_\varepsilon(\tilde \Psi_\alpha)\leq \alpha} \Beta_\varepsilon(\Theta_a(r_\varepsilon),\tilde \Psi_\alpha),
$$
for some radius $r_\varepsilon:=r_\varepsilon^\star >0$. Indeed,  if 
\begin{equation}
\inf_{\tilde \Psi_\alpha: \,\Alpha_\varepsilon(\tilde \Psi_\alpha)\leq \alpha} \Beta_\varepsilon(\Theta_a(r_\varepsilon^\star),\tilde \Psi_\alpha) \geq \beta,
\label{eq:Fig1-lb}
\end{equation}
for some $r_\varepsilon^\star >0$, then 
$$
\tilde r_{\varepsilon} \geq r_\varepsilon^\star.
$$

\noindent
{\bf Upper bound:} We first construct an $\alpha$-level test $\Psi_{\alpha}$. Then, we are looking for a radius $r_{\varepsilon}>0$ such that, uniformly over all $\theta \in \mathcal{E}_a$, 
$$
\|\theta\| > r_{\varepsilon} \quad \Rightarrow \quad \mathbb{P}_\theta(\Psi_\alpha=0) \leq \beta.
$$
It is then evident that
\begin{equation}
\Beta_\varepsilon(\Theta_a(r_\varepsilon),\Psi_\alpha)  \leq \beta \quad \mbox{implying that} \quad \ r_\varepsilon(\mathcal{E}_a,\Psi_\alpha,\beta) \leq r_{\varepsilon}.
\label{eq:Fig1-ub}
\end{equation}
Finally, if $r_{\varepsilon} \leq \mathcal{C} r_\varepsilon^\star$ for some $\mathcal{C} \geq1$, it then follows immediately that
$$
r_\varepsilon(\mathcal{E}_a,\Psi_\alpha,\beta) \leq \mathcal{C} r_\varepsilon^\star.
$$
(Note that the $\alpha$-level test $\Psi_{\alpha}$ constructed above is powerful according to Definition 
\ref{def:powerFanis}.)\\

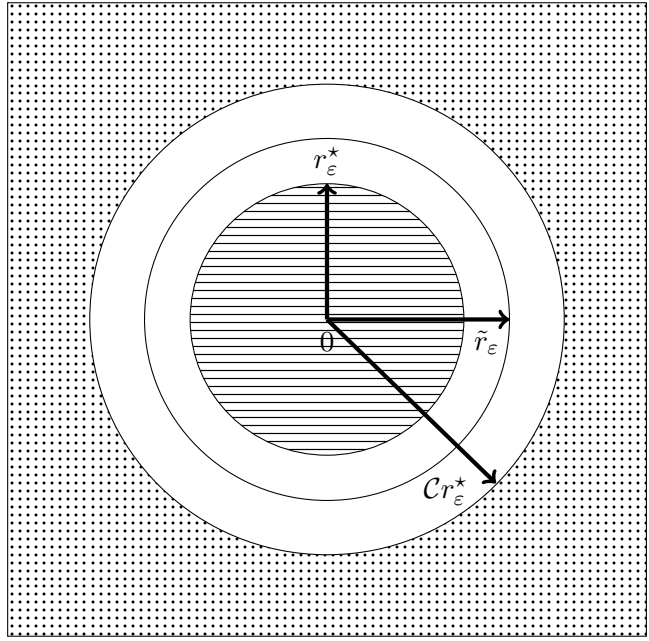
\begin{figure}
\begin{center}
\begin{tikzpicture}[scale=1.2]
\filldraw[pattern=dots,even odd rule](-3.5,-3.5) rectangle (3.5,3.5) (0,0) circle (2.6cm);
\draw[pattern=horizontal lines] (0,0) circle (1.5cm);
\draw (0,0) circle (2cm);
\draw (0,-0.25) node {0};
\draw [ultra thick,->] (0,0) -- (2,0);
\draw (1.75,-0.25) node {$\tilde r_\varepsilon$};
\draw [ultra thick, ->] (0,0) -- (0,1.5);
\draw (0,1.75) node {$r_\varepsilon^\star$};
\draw [ultra thick, ->] (0,0) -- (1.85,-1.81);
\draw (1.3, -1.9) node {$\mathcal{C} r_\varepsilon^\star$};
\end{tikzpicture}
\end{center}
\caption{\em According to the lower bound (\ref{eq:Fig1-lb}), minimax signal detection is not possible inside the circle with center $0$ and radius $r_\varepsilon^\star$. According to the upper bound (\ref{eq:Fig1-ub}), for any given $\alpha \in ]0,1[$, the maximal second kind error probability $\Beta_\varepsilon(\Theta_a(r_\varepsilon),\Psi_\alpha)$ of an $\alpha$-level test $\Psi_\alpha$ can be controlled by a prescribed level $\beta \in ]0,1-\alpha[$ outside the circle with center $0$ and radius $\mathcal{C} r_\varepsilon^\star$, for some $\mathcal{C} \geq 1$. The `optimal frontier' is determined by the circle with center $0$ and radius $\tilde r_\varepsilon$ $($i.e., the minimax separation radius$)$.}
\label{fig:Fig1}
\end{figure}

Figure \ref{fig:Fig1} illustrates the areas where, according to Definitions \ref{eq:minimax_separation_radius}--\ref{def:powerFanis}, minimax signal detection can, or cannot, be possible.\\

We stress at this point that the quantity 
$$
\inf_{\tilde \Psi_\alpha: \,\Alpha_\varepsilon(\tilde \Psi_\alpha)\leq \alpha} \Beta_\varepsilon(\Theta_a(r_\varepsilon),\tilde \Psi_\alpha)
$$
needed to bound from below in the above discussion is precisely the minimax second kind error probability to be introduced in the asymptotic approach that we elaborate in the following section.

\subsection{The asymptotic approach}
\label{asympt_general}
Let $\alpha \in ]0,1[$ be fixed and let $r_\varepsilon >0$ be a given radius. 

\begin{definition}
\label{mskep-F}
The minimax second kind error probability is defined as
$$ \Beta_ {\varepsilon,\alpha}(\Theta_a(r_\varepsilon)) := \inf_{\tilde \Psi_\alpha: \,\Alpha_\varepsilon(\tilde \Psi_\alpha)\leq \alpha} \; \Beta_\varepsilon(\Theta_a(r_\varepsilon),\tilde \Psi_\alpha).$$ 
\end{definition}

Given a radius $r_\varepsilon >0$, the minimax second kind error probability $\Beta_ {\varepsilon,\alpha}(\Theta_a(r_\varepsilon))$ characterizes the minimax testing performances over all $\alpha$-level tests $\tilde \Psi_\alpha$ for signal detection problem (\ref{testing_pb1}). In other words, it corresponds to the lowest maximal second kind error probability over the set $\Theta_a(r_\varepsilon)$. In particular, one would like to identify the different possible values of the radius $r_\varepsilon>0$ such that the minimax second kind error probability $\Beta_{\varepsilon,\alpha}(\Theta_a(r_\varepsilon))$ tends to $0$ or to a constant or to $1$, as $\varepsilon$ tends to $0$.

\begin{definition}
\label{def:separation}
The term $\bar r_\varepsilon:=\bar r_\varepsilon(\mathcal{E}_a,\alpha) >0$ is called the minimax separation rate if, for any given $r_{\varepsilon}>0$,
$$ \Beta_ {\varepsilon,\alpha}(\Theta_a(r_\varepsilon)) =1-\alpha +o_\varepsilon(1) \quad \mathrm{if} \quad \frac{r_\varepsilon}{\bar r_\varepsilon} \rightarrow 0 \quad \mathrm{as} \quad \varepsilon\rightarrow 0,
$$
and 
$$\Beta_ {\varepsilon,\alpha}(\Theta_a(r_\varepsilon)) = o_\varepsilon(1) \quad \mathrm{if} \quad \frac{r_\varepsilon}{\bar r_\varepsilon} \rightarrow +\infty \quad \mathrm{as} \quad \varepsilon\rightarrow 0.
$$
\end{definition}

The minimax separation rate $\bar r_\varepsilon$ identifies, in some sense, the frontiers between detectable and undetectable signals. In other words, it means that, for small $\varepsilon >0$, one can detect {\em all}\, $\theta \in \Theta_a(r_\varepsilon)$ for which the ratio $r_\varepsilon/\bar r_\varepsilon$ is large. On the other hand, if, for small $\varepsilon >0$, the ratio $r_\varepsilon/\bar r_\varepsilon$ is small, it is then {\em impossible} to distinguish $H_0$ from $H_1$ with small maximal second kind error probability $\Beta_{\varepsilon,\alpha}(\Theta_a(r_\varepsilon))$. \\

In practice, given an $\alpha$-level test $\Psi_\alpha$, it might be useful, for small $\varepsilon >0$, to compare its maximal second kind error probability $\Beta_\varepsilon(\Theta_a(r_\varepsilon),\Psi_\alpha)$ to the minimax second kind error probability $\Beta_{\varepsilon,\alpha}(\Theta_a(r_\varepsilon))$. Hence, the following definition is appropriate.

\begin{definition}
\label{def:am}
An $\alpha$-level test $\Psi_\alpha$ is said to be \\

(i) asymptotical minimax consistent if, for any given $r_\varepsilon >0$,
$$
\Beta_\varepsilon(\Theta_a(r_\varepsilon),\Psi_\alpha) = o_\varepsilon(1) \quad \mathrm{if} \quad \frac{r_\varepsilon}{\bar r_\varepsilon} \rightarrow +\infty \quad \mathrm{as} \quad \varepsilon\rightarrow 0.
$$

(ii) asymptotical minimax if, for any given $r_\varepsilon >0$, 
$$\Beta_\varepsilon(\Theta_a(r_\varepsilon),\Psi_\alpha) = \Beta_ {\varepsilon,\alpha}(\Theta_a(r_\varepsilon)) + o_\varepsilon(1).$$
\end{definition}

Regarding Definition \ref{def:am}, given an $\alpha$-level test $\Psi_\alpha$, item (i) provides a weak condition in the sense that, for small $\varepsilon >0$, one can detect all $\theta \in \Theta_a(r_\varepsilon)$ for which the ratio $r_\varepsilon/\bar r_\varepsilon$ is large. On the other hand, item (ii) refers to a strong condition in the sense that one needs to asymptotically attain the minimax second kind error probability $\Beta_{\varepsilon,\alpha}(\Theta_a(r_\varepsilon))$.\\
 
In this setting, the point of view is asymptotic. The performance of any testing procedure is investigated as $\varepsilon$ tends to 0.  Nevertheless, such a point of view allows, sometimes, to provide a precise description of the asymptotic value for the minimax separation rate $\bar r_\varepsilon$. In particular, one can, in some cases, determine sharp asymptotics of Gaussian type for the minimax second kind error probability $\Beta_{\varepsilon,\alpha}(\Theta_a(r_\varepsilon))$. 

\begin{definition} 
\label{sharp_asy}
The minimax second kind error probability $\Beta_ {\varepsilon,\alpha}(\Theta_a(r_\varepsilon))$ is said to possess a sharp asymptotic of Gaussian type if it has an asymptotic Gaussian shape, i.e., if there exists a function $\nu(r_{\varepsilon}) \in \, ]-\infty, \Phi^{-1}(1-\alpha)]$ (that should be determined later on) such that
$$ \Beta_ {\varepsilon,\alpha}(\Theta_a(r_\varepsilon)) = \Phi(\nu(r_{\varepsilon})) + o_\varepsilon(1),$$ 
where $\Phi$ denotes the distribution function of the standard Gaussian distribution.
\end{definition}

Sharp asymptotics of Gaussian type for the minimax second kind error probability $\Beta_ {\varepsilon,\alpha}(\Theta_a(r_\varepsilon))$ have been observed in particular settings (see e.g., \cite{IS_2003} and references therein).\\

We present below a general strategy for obtaining the minimax separation rate $\bar r_\varepsilon$ and sharp asymptotics of Gaussian type for the minimax second kind error probability $\Beta_ {\varepsilon,\alpha}(\Theta_a(r_\varepsilon))$. Given an ellipsoid $\mathcal{E}_a$, this amounts to investigate the construction of both lower and upper bounds on $\Beta_{\varepsilon,\alpha}(\Theta_a(r_\varepsilon))$.\\

\noindent
{\bf Lower bound:} Find a radius $r_{\varepsilon,1}>0$  such that, for any given $r_\varepsilon >0$,
$$\Beta_ {\varepsilon,\alpha}(\Theta_a(r_\varepsilon)) \geq 1-\alpha +o_{\varepsilon}(1) \quad \mathrm{if} \quad \frac{r_\varepsilon}{r_{\varepsilon,1}} \rightarrow 0 \quad \mathrm{as} \quad \varepsilon\rightarrow 0.
$$
If possible, one may also want to determine the shape of $\Beta_{\varepsilon,\alpha}(\Theta_a(r_\varepsilon))$, i.e., to find a function $\nu_1(r_{\varepsilon}) \in \, ]-\infty, \Phi^{-1}(1-\alpha)]$ such that, for any given $r_\varepsilon >0$,
$$\Beta_ {\varepsilon,\alpha}(\Theta_a(r_\varepsilon)) \geq \Phi(\nu_1(r_{\varepsilon})) +o_{\varepsilon}(1).
$$

\noindent
{\bf Upper bound:} Given an $\alpha$-level test $\Psi_{\alpha}$, find a radius $r_{\varepsilon,2}>0$ such that, for any given $r_\varepsilon >0$,
$$
\Beta_\varepsilon(\Theta_a(r_\varepsilon),\Psi_\alpha) = o_\varepsilon(1) \quad \mathrm{if} \quad \frac{r_\varepsilon}{r_{\varepsilon,2}} \rightarrow +\infty \quad \mathrm{as} \quad \varepsilon\rightarrow 0.
$$
Additionally, one may again want to determine the shape of $\Beta_\varepsilon(\Theta_a(r_\varepsilon),\Psi_\alpha)$,  i.e., to find a function $\nu_2(r_{\varepsilon}) \in \, ]-\infty, \Phi^{-1}(1-\alpha)]$ such that, for any given $r_\varepsilon >0$,
$$
\Beta_\varepsilon(\Theta_a(r_\varepsilon),\Psi_\alpha) \leq \Phi(\nu_2(r_{\varepsilon})) +o_{\varepsilon}(1).
$$

If the $\alpha$-level test $\Psi_{\alpha}$ is such that $r_{\varepsilon,1}/r_{\varepsilon,2} =\mathcal{O}_\varepsilon(1)$, then, obviously, $\bar{r}_{\varepsilon}/r_{\varepsilon,1} =\mathcal{O}_\varepsilon(1)$. It means that, according to Definition \ref{def:separation}, either $r_{\varepsilon,1}$ or $r_{\varepsilon,2}$ correspond to the minimax separation rate $\bar r_\varepsilon$. Furthermore, in the case when $\nu_1(r_{\varepsilon})/\nu_2(r_{\varepsilon})=1+o_{\varepsilon}(1)$, then, according to Definition \ref{sharp_asy}, we get sharp asymptotics of Gaussian type for the minimax second kind error probability $\Beta_ {\varepsilon,\alpha}(\Theta_a(r_\varepsilon))$, with $\nu(\cdot)=\nu_1(\cdot)$. \\

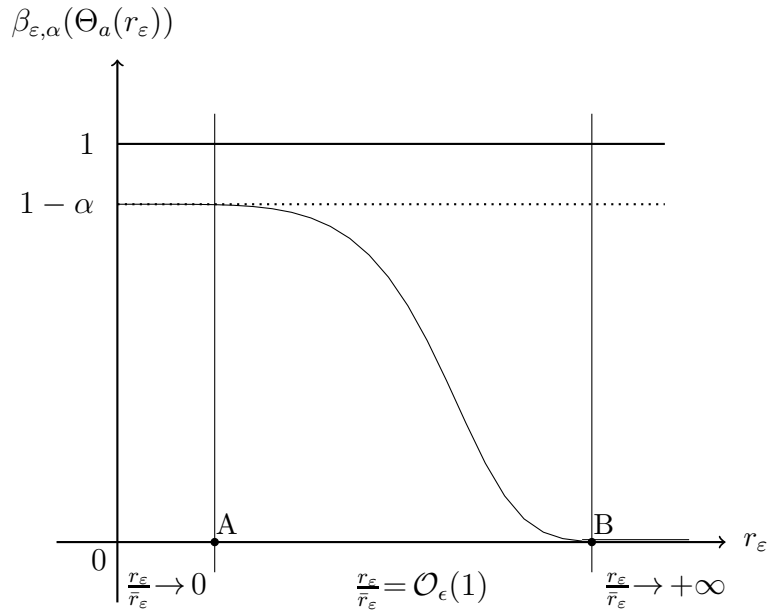
\begin{figure}
\begin{center}
\begin{tikzpicture}[scale=1.60]
\draw [thick,->] (-0.5,0) -- (5,0);
\draw (5.25,0) node {$r_\varepsilon$};
\draw [thick, ->] (0,-0.5) -- (0,4);
\draw (-0.15,-0.15) node {\large $0$};
\draw (-0.2,4.3) node {{\large $\beta_{\varepsilon,\alpha}(\Theta_a(r_\varepsilon))$}};
\draw [thick] (-0,3.3) -- (4.5,3.3);
\draw (-0.25,3.3) node {\large $1$};
\draw (-0.5,2.8) node {\large $1-\alpha$};
\draw [thick,dotted] (-0,2.8) -- (4.5,2.8);
\draw [black,domain=0:3.82,thin] plot (\x, {2*2.8*exp(-0.008*\x*\x*\x*\x*\x)/(1+exp(-0.00*\x-0.008*\x*\x*\x*\x*\x))});
\draw [black,domain=3.82:4.7,thin] plot (\x, {0.02});
\draw [thin] (0.8,-0.25) -- (0.8,3.55);
\draw (0.9,0.15) node {A};
\draw[fill] (0.8,0) circle (0.03);
\draw [thin] (3.9,-0.25) -- (3.9,3.55); 
\draw (4,0.15) node {B};
\draw[fill] (3.9,0) circle (0.03);
\draw (0.4,-0.4) node {{\large $\frac{r_\varepsilon}{\bar r_\varepsilon} \hspace{-0.1cm}\rightarrow \hspace{-0.05cm} 0$}};
\draw (4.5,-0.4) node {{\large $\frac{r_\varepsilon}{\bar r_\varepsilon} \hspace{-0.1cm}\rightarrow \hspace{-0.05cm} +\infty$}};
\draw (2.5,-0.4) node {{\large $\frac{r_\varepsilon}{\bar r_\varepsilon} \hspace{-0.1cm}= \hspace{-0.05cm} \mathcal{O}_\epsilon(1)$}};
\end{tikzpicture}
\end{center}
\caption{\textit{The interval $[0,A]$ $($resp. $[B,+\infty[$$)$ illustrates the area where $\frac{r_\varepsilon}{\bar r_\varepsilon} \rightarrow 0$ $($resp. $\frac{r_\varepsilon}{\bar r_\varepsilon} \rightarrow +\infty$$)$ as $\varepsilon \rightarrow 0$, i.e., where the minimax second kind error probability $\beta_{\varepsilon,\alpha}(\Theta_a(r_\varepsilon))$ satisfies $\beta_{\varepsilon,\alpha}(\Theta_a(r_\varepsilon))= 1-\alpha+ o_\epsilon(1)$ $($resp. $\beta_{\varepsilon,\alpha}(\Theta_a(r_\varepsilon))=o_\epsilon(1)$$)$ $($see Definition \ref{def:separation}$)$. The interval $[A,B]$ determines the frontier for the minimax separation rate $\bar r_\varepsilon$. In particular, inside this area, sharp asymptotics of Gaussian type (solid curve) for $\beta_{\varepsilon,\alpha}(\Theta_a(r_\varepsilon))$ are feasible, i.e., there exists a function $\nu(r_{\varepsilon}) \in \, ]-\infty, \Phi^{-1}(1-\alpha)]$ such that $\beta_{\varepsilon,\alpha}(\Theta_a(r_\varepsilon)) = \Phi(\nu(r_\varepsilon))+o_\varepsilon(1)$ $($see Definition \ref{sharp_asy}$)$.}}
\label{fig:Fig2}
\end{figure}

Figure \ref{fig:Fig2} illustrates the areas where, according to Definition \ref{def:separation}, minimax signal detection can, or cannot, be possible. It also illustrates, according to Definition \ref{sharp_asy}, the area where sharp asymptotics of Gaussian type for the minimax second kind error probability $\Beta_ {\varepsilon,\alpha}(\Theta_a(r_\varepsilon))$ are feasible.

\subsection{A brief motivation}

Although the considered minimax signal detection problem (\ref{testing_pb1}) is the same for both approaches (non-asymptotic and asymptotic), the way the optimality of the considered testing procedures is measured differs. \\

In the non-asymptotic setting, the statistician sets in advance some prescribed values $\alpha, \beta \in ]0,1[$. Then, the goal is to find `optimal' (non-asymptotic) separation conditions for $H_0$ and $H_1$ that allow a precise (non-asymptotic) control of the first kind error probability and maximal second kind error probability by $\alpha$ and $\beta$, respectively. On the other hand, in the asymptotic setting, the aim is slightly different. Given any $r_\varepsilon >0$, the goal is to measure the best possible associated maximal second kind error probability of an (asymptotical) $\alpha$-level test and to (asymptotically) determine whether it tends to $1-\alpha$ or to $0$, as the noise level $\varepsilon$ tends to 0.\\

In order to study the signal detection problem (\ref{testing_pb1}), from a minimax point of view, different testing methodologies have been developed over the years that strongly depend on the two considered signal detection paradigms. We refer to, e.g., \cite{Baraud}, \cite{LLM_2011}, \cite{LLM_2012}, for the non-asymptotic paradigm, and to, e.g., \cite{IngsterI}, \cite{IngsterII}, \cite{IngsterIII}, \cite{IS_2003}, \cite{ISS_2012}, for the asymptotic paradigm. Unsurprisingly, the results in these studies are coherent (i.e., the associated minimax separation radii $\tilde{r}_\varepsilon$ and minimax separation rates $\bar{r}_\varepsilon$ are {\em asymptotically} equivalent, as $\varepsilon$ tends to 0). Indeed, one can formally prove (using the respective definitions) that $\tilde r_\epsilon / \bar r_\epsilon = O_\epsilon(1)$ as $\epsilon\rightarrow 0$. \\

In the sequel, we propose below a unified treatment for the study of the minimax separation radius $\tilde r_\varepsilon$ (non-asymptotic approach) and the minimax separation rate $\bar r_\varepsilon$ (asymptotic approach). We compare the construction of their lower and upper bounds and point out similarities in both settings (Sections  \ref{s:control} and \ref{connection}). In particular, tools constructed in the non-asymptotic paradigm can be used in order to draw  conclusions in the asymptotic paradigm and vice-versa. 
%we prove that the minimax separation radius $\tilde r_\varepsilon$ and the minimax separation rate $\bar r_\varepsilon$ are, under mild conditions,  asymptotically equivalent (Section \ref{subsec:explicit_seq}). 
In other words, one can perform asymptotic analysis for non-asymptotic testing procedures and investigate non-asymptotic performances for asymptotic testing procedures. This will be demonstrated later on, when explicit sequences $a=(a_j)_{j\in \mathbb{N}}$ and $b=(b_j)_{j\in \mathbb{N}}$ are at hand (see Section \ref{s:example}).

\section{Control of the lower and upper bounds}
\label{s:control}

\subsection{Control of the lower bounds}
\label{s:lower_bound}
One of the main issues of minimax signal detection is to establish lower bounds for the minimax separation radius $\tilde r_\varepsilon$ (non-asymptotic approach) and the minimax separation rate $\bar r_\varepsilon$ (asymptotic approach). In both approaches, this amounts to determine the values of the available signal for which $H_0$ and $H_1$ cannot be separated with prescribed minimax second kind error probability $\Beta_{\varepsilon,\alpha}(\Theta_a(r_\varepsilon))$. \\

More formally, we are interested to bound from below the minimax second kind error probability $\Beta_{\varepsilon,\alpha}(\Theta_a(r_\varepsilon))$.
In particular, an interesting question is to investigate the smallest possible value of the radius $r_\varepsilon >0$ for which $\Beta_{\varepsilon,\alpha}(\Theta_a(r_\varepsilon))$ can be, following the non-asymptotic or asymptotic approaches, (asymptotically) lower bounded by $\beta \in \, ]0,1-\alpha[$ or tends to $1-\alpha$.\\

A possible way to achieve this goal is to consider a (prior) probability measure $\pi$ on the set associated with $H_1$, i.e., a probability measure $\pi$ on the set  $\Theta_a(r_\varepsilon)$ (see, e.g., \cite{Baraud}, \cite{IS_2003}). Then, it is easily verified that
\begin{eqnarray*}
\Beta_{\varepsilon,\alpha}(\Theta_a(r_\varepsilon))
& \geq & \inf_{\tilde \Psi_\alpha: \,\Alpha_\varepsilon(\tilde \Psi_\alpha)\leq \alpha} \mathbb{P}_\pi (\tilde \Psi_\alpha =0)\\ 
& = & \inf_{\tilde \Psi_\alpha: \,\Alpha_\varepsilon(\tilde \Psi_\alpha)\leq \alpha} \left[ \mathbb{P}_0(\tilde \Psi_\alpha =0) + \mathbb{P}_\pi (\tilde \Psi_\alpha =0) -\mathbb{P}_0(\tilde \Psi_\alpha =0)\right]\\
& \geq &  \inf_{\tilde \Psi_\alpha: \,\Alpha_\varepsilon(\tilde \Psi_\alpha)\leq \alpha} \left[1-\alpha - \left| \mathbb{P}_\pi (\tilde \Psi_\alpha =0) -\mathbb{P}_0(\tilde \Psi_\alpha =0)\right|  \right]  \\
& \geq & 1- \alpha - \sup_{A:  \; \mathbb{P}_0(A)\leq \alpha} |\mathbb{P}_\pi(A)-\mathbb{P}_0(A) |\\
& \geq & 1- \alpha - \sup_{A  \in {\cal A}} |\mathbb{P}_\pi(A)-\mathbb{P}_0(A) |\\
& = & 1-\alpha - V(\mathbb{P}_\pi,\mathbb{P}_0),
\end{eqnarray*}
where $$V(\mathbb{P}_\pi,\mathbb{P}_0):=\sup_{A \in {\cal A}}|\mathbb{P}_\pi(A)-\mathbb{P}_0(A) |$$ denotes the total variation norm between the two probability measures $\mathbb{P}_0$ and $\mathbb{P}_\pi=\int \mathbb{P}_{\theta}\,d\pi(\theta)$), and  ${\cal A}$ denotes the $\sigma$-field of the underlying probability space. Assuming that $\mathbb{P}_\pi$ is absolutely continuous with respect to $\mathbb{P}_0$, using first the Scheff\'e Theorem (see, e.g., \cite{Tsybakov_2009}, Lemma 2.1) and then the Cauchy-Schwarz inequality, it can be seen that
\begin{eqnarray*}
V(\mathbb{P}_\pi,\mathbb{P}_0) & := & \sup_{A \in {\cal A}} |\mathbb{P}_\pi(A)-\mathbb{P}_0(A) |\\
& = & \frac{1}{2} \int \left| \frac{d\mathbb{P}_\pi}{dy}(y) - \frac{d\mathbb{P}_0}{dy}(y)   \right| dy\\
& = & \frac{1}{2} \int \left| \frac{d\mathbb{P}_\pi}{d\mathbb{P}_0}(y) - 1   \right| d\mathbb{P}_0(y)\\
& \leq & \frac{1}{2} \left( \mathbb{E}_{0} ( | L_\pi(Y)-1|^2) \right)^{1/2},
\end{eqnarray*}
where $L_\pi(Y)$ denotes the likelihood ratio between the two measures $\mathbb{P}_\pi$ and $\mathbb{P}_0$, and $\mathbb{E}_0$ denotes the expectation with respect to $\mathbb{P}_0$.  Combining the above arguments, we obtain the following lower bound
\begin{equation}
\Beta_{\varepsilon,\alpha}(\Theta_a(r_\varepsilon)) \geq 1- \alpha -\frac{1}{2} \left( \mathbb{E}_0[ L_{\pi}^2(Y)] -1 \right)^{1/2}.
\label{eq:step1}
\end{equation}

The construction of the lower bound for the minimax second kind error probability $\Beta_{\varepsilon,\alpha}(\Theta_a(r_\varepsilon))$ developed in (\ref{eq:step1}) heavily relies on the construction of a prior $\pi$ on the set $\Theta_a(r_\varepsilon)$. Given some sequence $\theta=(\theta_j)_{j\in\mathbb{N}}$, which will be made explicit below, we consider the symmetric prior $\pi$ defined as
\begin{equation}
\pi = \prod_{j\in \mathbb{N}} \pi_j \quad \mathrm{with} \quad \pi_j=\frac{1}{2} (\delta_{-\theta_j} + \delta_{\theta_j}) \quad \forall \, j\in \mathbb{N}.
\label{eq:prior}
\end{equation}
(Note that $\pi(\Theta_a(r_\varepsilon)=1$.)
Since the $\xi_j$ are standard Gaussian random variables, we get, after some technical algebra (see \cite{Baraud} p. 596 or \cite{ISS_2012}, supplementary material, Section 11.1), that
\begin{equation}
\mathbb{E}_0[ L_{\pi}^2(Y)] = \prod_{j\in\mathbb{N}} \cosh (b_j^2\theta_j^2/\varepsilon^2) \leq \exp \left( \frac{1}{2\varepsilon^4} \sum_{j\in \mathbb{N}} b_j^{4}\theta_j^{4} \right) := \exp(u_{\varepsilon,\theta}^2).
\label{eq:control_LR}
\end{equation}

It is worth pointing out that the construction of the lower bound for the minimax separation radius $\tilde r_\varepsilon$ (non-asymptotic approach) and the minimax separation rate $\bar r_\varepsilon$(asymptotic approach) are then both related to the study of either $\mathbb{E}_0[ L_{\pi}^2(Y)]$ or its corresponding upper bound (\ref{eq:control_LR}).\\

Two different interesting regimes at this point can be immediately deduced:
\begin{itemize}
\item First, $\mathbb{E}_0[ L_{\pi}^2(Y)]$ tends to $1$ as $\varepsilon \rightarrow 0$. In such a case, the minimax second kind error probability $\Beta_{\varepsilon,\alpha}(\Theta_a(r_\varepsilon))$ is asymptotically lower bounded by $1-\alpha$, i.e., $\Beta_{\varepsilon,\alpha}(\Theta_a(r_\varepsilon)) \geq 1-\alpha +o_\varepsilon(1)$. 
\item Second, $\mathbb{E}_0[ L_{\pi}^2(Y)]$ can be upper bounded by a constant. In this case, the minimax second kind error probability $\Beta_{\varepsilon,\alpha}(\Theta_a(r_\varepsilon)) $ is also lower bounded by a constant. An interesting situation corresponds to the case where $\mathbb{E}_0[ L_{\pi}^2(Y)] \leq 1+4(1-\alpha-\beta)^2$ for some $\beta\in ]0,1-\alpha[$. Then, $\Beta_{\varepsilon,\alpha}(\Theta_a(r_\varepsilon)) \geq \beta$. 
\end{itemize}

Moreover, a more delicate study of the term $u^2_{\varepsilon,\theta}:=\frac{1}{2\varepsilon^4} \sum_{j\in \mathbb{N}} b_j^{4}\theta_j^{4}$ in (\ref{eq:control_LR}) allows, under certain conditions to be made precise later on, to study sharp asymptotics of Gaussian type for the minimax second kind error probability $\Beta_ {\varepsilon,\alpha}(\Theta_a(r_\varepsilon))$.\\

We discuss below the two different strategies that have been investigated in the literature.

\subsubsection{Non-asymptotic control}
\label{s:lb_nasymp}
The following control has been proposed by \cite{Baraud} in the direct setting and it has been generalized to the inverse setting by \cite{LLM_2012}. The main idea
consists of finding an explicit sequence $\theta^0=(\theta_j^0)_{j\in \mathbb{N}}$ and a radius $r_\varepsilon >0$ which satisfy the following three requirements:
\begin{itemize}
\item $\| \theta^0 \| \geq r_\varepsilon$,
\item $\exp \left( \frac{1}{2\varepsilon^4} \sum_{j\in \mathbb{N}} b_j^{4}(\theta_j^0)^{4} \right)=1+4(1-\alpha-\beta)^2$,
\item $\theta^0 \in \mathcal{E}_a$.
\end{itemize}
To this end, one can consider, for instance, the sequence $\theta^0$ defined as
\begin{equation} 
\theta_j^0 := \frac{r_\varepsilon \varepsilon^2 b_j^{-2}}{\left( \varepsilon^4 \sum_{k=1}^D b_k^{-4} \right)^{1/2}} \ \ \forall \, j\in \left\lbrace 1,\dots, D \right\rbrace \quad \mathrm{and} \quad \theta_j^0 =0 \ \ \forall \, j>D, 
\label{eq:step2}
\end{equation}
for some (finite) parameter $D \in \mathbb{N}$ (called the bandwidth), that possibly depends on $\varepsilon>0$. \\

It is evident that $\|\theta^0\| = r_\varepsilon$. Furthermore, taking into account (\ref{eq:control_LR}), we get 
\begin{equation} 
\mathbb{E}_0[ L_{\pi}^2(Y)]  \leq \exp \left( \frac{1}{2\varepsilon^4} \sum_{j\in \mathbb{N}} b_j^{4}(\theta_j^0)^{4} \right) = \exp\left[ \frac{r_\varepsilon^4}{2\varepsilon^4 \sum_{j=1}^D b_j^{-4}} \right] = 1+4(1-\alpha-\beta)^2,
\label{eq:step3}
\end{equation}
as soon as
$$
r_\varepsilon^2 = r_{\varepsilon,D}^2 := c(\alpha,\beta) \varepsilon^2 \sqrt{ \sum_{j=1}^D b_j^{-4}},
$$
where 
\begin{equation}
c(\alpha,\beta)= (2 \ln(1+4(1-\alpha-\beta)^2))^{1/4}>0.
\label{clement_c}
\end{equation}
In order to conclude, it remains to choose an appropriate $D \in \mathbb{N}$ such that $\theta^0 \in \mathcal{E}_a$. To this end, note that
$$ \sum_{j\in \mathbb{N}} a_j^2 (\theta_j^0)^2 \leq a_D^{2} \sum_{j=1}^D  ( \theta_j^0)^2 = a_D^2 r_{\varepsilon,D}^2 \leq 1 \quad \text{as soon as} \quad \ r_{\varepsilon,D}^2 \leq a_D^{-2} . $$
Hence, if we define
\begin{equation}
r_{\varepsilon,\star}^2 := \sup_{D\in \mathbb{N}} \left[ c(\alpha,\beta) \varepsilon^2 \sqrt{\sum_{j=1}^D b_j^{-4}} \wedge a_D^{-2} \right],
\label{bi_nasymp1}
\end{equation}
we get, using (\ref{eq:step1})-(\ref{eq:step3}), 
$$ \Beta_{\varepsilon,\alpha}(\Theta_a(r_{\varepsilon,\star})) \geq \beta,$$
which means that the minimax separation radius $\tilde r_{\varepsilon}$ satisfies
$$\tilde r_{\varepsilon} \geq r_{\varepsilon,\star}.$$

This corresponds to a non-asymptotic lower bound on the minimax separation radius $\tilde r_{\varepsilon}$. The main advantage of such a bound is that it provides a precise description of the area where minimax signal detection is impossible with prescribed values $\alpha, \beta \in ]0,1[$. 

\subsubsection{Asymptotic control}
\label{subsubsec:ac-final-1}

In the previous (non-asymptotic) approach, the main idea was to construct an explicit sequence $\theta$ and to control $\mathbb{E}_0[ L_{\pi}^2(Y)]$. In the asymptotic approach, one instead starts from (\ref{eq:control_LR}) and find the smallest possible value of $u_{\varepsilon,\theta}^2$ for which $\theta \in \Theta_a(r_\varepsilon)$. In other words, the idea is to choose a sequence $\bar \theta:=\bar \theta(r_\varepsilon)$ as the solution of the following extremal problem
\begin{equation}
\bar \theta(r_\varepsilon) := {\arginf}_{\theta \in \Theta_a(r_\varepsilon)} \left\lbrace u_{\varepsilon,\theta}^2:= \frac{1}{2\varepsilon^4} \sum_{k\in \mathbb{N}} b_k^{4}\theta_k^4\right\rbrace. 
\label{extrem_pb}
\end{equation}
In the following, we will denote the solution of the extremal problem (\ref{extrem_pb}) as
\begin{equation}
u_\varepsilon(r_\varepsilon) := u_{\varepsilon,\bar \theta(r_\varepsilon)}:={\inf}_{\theta \in \Theta_a(r_\varepsilon)} \left\lbrace  \frac{1}{2\varepsilon^4} \sum_{k\in \mathbb{N}} b_k^{4}\theta_k^4\right\rbrace.
\label{extrem_sol}
\end{equation}

This idea has been in particular developed in the series of papers \cite{IngsterI}, \cite{IngsterII}, \cite{IngsterIII}, or, more recently, in \cite{ISS_2012}, in an inverse problem framework. The cases of interest correspond to the setting where $u_\varepsilon^2(r_\varepsilon)$ either tends to zero or is bounded by a constant. In that case, one can find the solution of $u_\varepsilon(r_\varepsilon)$ in (\ref{extrem_sol}) using, for instance, the standard methodology of Lagrange multipliers.   \\

Let $\alpha \in ]0,1[$ be fixed. We immediately see from (\ref{eq:step1}) and (\ref{eq:control_LR}) that if $u_\varepsilon(r_\varepsilon)=o_\varepsilon(1)$, then 
\begin{equation}
\label{eq:lower_bound_zero}
\Beta_{\varepsilon,\alpha}(\Theta_a(r_\varepsilon)) \geq 1-\alpha +o_\varepsilon(1).
\end{equation} 
\vspace{-0.4cm}

The interesting situation, however, arises when  $u_\varepsilon(r_\varepsilon) = \mathcal{O}_\varepsilon(1)$. It allows a more accurate study to asymptotically precise the shape of the minimax second kind error probability $\Beta_{\varepsilon,\alpha}(\Theta_a(r_\varepsilon))$. In particular, if $u_\varepsilon(r_\varepsilon) = \mathcal{O}_\varepsilon(1)$, it can be established that
\begin{equation}
\label{fanis:loglikeasymp}
 \ln(L_\pi(Y)) = -\frac{u_\varepsilon^2(r_\varepsilon)}{2} + u_\varepsilon(r_\varepsilon) \xi_\varepsilon + \zeta_\varepsilon,
 \end{equation}
where $\xi_\varepsilon \rightarrow \xi \sim \mathcal{N}(0,1)$ and $\zeta_\varepsilon \rightarrow 0$ (in $\mathbb{P}_0$-probability) as $\varepsilon \rightarrow 0$ distribution (see Section 4.3.1 of \cite{IS_2003} or the proof of Theorem 4.1 of \cite{ISS_2012}, supplementary material, Section 11.1). By a standard change of probability measure, it follows that
$$
\Beta_{\varepsilon,\alpha}(\Theta_a(r_\varepsilon)) \geq \mathbb{E}_{\pi} (1-\psi^\star)=
\mathbb{E}_0 (\exp(\ln(L_{\pi}(Y))) (1-\psi^\star)),
$$
where $\psi^\star$ is the likelihood ratio test defined as $\psi^\star= \mathbf{1}_{\lbrace \ln(L_{\pi}(Y) )> t^\star_{1-\alpha}\rbrace}$ with $t^\star_{1-\alpha}$ being the $(1-\alpha)$-quantile of the distribution of $\ln(L_{\pi}(Y) )$ under $H_0$. Hence, in view of (\ref{fanis:loglikeasymp}), it is easily seen that
$$
t^\star_{1-\alpha}= -\frac{u_\varepsilon^2(r_\varepsilon)}{2} + u_\varepsilon(r_\varepsilon) t_{1-\alpha} +o_\varepsilon(1),
$$
where $t_{1-\alpha}$ refers to the $(1-\alpha)$-quantile of the standard Gaussian distribution. Moreover, using the mean value theorem, it follows that
\begin{eqnarray}
\Beta_{\varepsilon,\alpha}(\Theta_a(r_\varepsilon)) \geq \mathbb{E}_{\pi} (1-\psi^\star) &=& 
\mathbb{E}_0 (\exp(\ln(L_{\pi}(Y))) (1-\psi^\star)) \nonumber \\ &=& 
\Phi(t_{1-\alpha} - u_\varepsilon(r_\varepsilon)) +o_\varepsilon(1).
\label{sharp_asymp}
\end{eqnarray}
%Therefore, according to Definition \ref{sharp_asy}, sharp asymptotics of Gaussian type for the lower bound of the minimax second kind error probability $\Beta_ {\varepsilon,\alpha}(\Theta_a(r_\varepsilon))$ is also obtained.\\
(Note that, in the particular case where $r_\varepsilon >0$ satisfies
$u_\varepsilon(r_\varepsilon)=t_{1-\alpha} - t_\beta$, 
then, it is immediately seen  that 
$\Beta_{\varepsilon,\alpha}(\Theta_a(r_\varepsilon))\geq \beta + o_\varepsilon(1)$.)

\begin{remark}
{\rm It is worth mentioning that one cannot determine at this point  the radius $r_{\varepsilon,1}>0$ (considered in the general strategy of Section \ref{asympt_general} for constructing lower bounds), unless the sequences $a=(a_j)_{j\in \mathbb{N}}$ and $b=(b_j)_{j\in\mathbb{N}}$ are explicitly given. We refer to, e.g., \cite{ISS_2012} for more details or to Section \ref{s:example} where a mildly ill-posed inverse problem is treated for illustrative purposes.}
\end{remark}
 
In the following section, we investigate upper bounds on the minimax separation radius $\tilde r_\varepsilon$ (non-asymptotic approach) and upper bounds on the minimax separation rate $\bar r_\varepsilon$ (asymptotic approach). In the latter setting, we also provide, under mild conditions, sharp asymptotics of Gaussian type for the minimax second kind error probability $\Beta_{\varepsilon,\alpha}(\Theta_a(r_\varepsilon))$. 

\subsection{Control of the upper bounds}

\subsubsection{A general testing methodology}
\label{s:general_test}
In this section, we construct appropriate tests and investigate the associated separation radius (non-asymptotic approach) and the maximal second kind error probability (asymptotic approach). Starting from signal detection problem (\ref{testing_pb1}), the underlying question is to decide whether we observe a signal or not. To this end, a possible approach is to construct an estimator $\hat d$ of $d:=\|\theta\|^2$  or $d:=\| b\theta\|^2$.  Indeed, the assertions $``\theta=0"$ and $``b\theta=0"$ are equivalent since the sequence $b=\{b_j\}_{j \in \mathbb{N}}$ is assumed strictly positive. We refer to \cite{LLM_2011} for an extended discussion on that subject. Then, one can use the following decision rule:
\begin{itemize}
\item if $\hat d$ is large enough (larger than a prescribed threshold which should be precisely quantified), we reject $H_0$,
\item If $\hat d$ is smaller than this threshold, we do not reject $H_0$.
\end{itemize}
In order to estimate $\| \theta \|^2$ (resp. $\| b\theta\|^2$), one can first construct a preliminary estimator of $\theta$ (resp. $b\theta$) and then take its squared norm. This idea has been widely investigated. We point out that, in general, the preliminary estimators cannot be directly plugged in order to estimate $\| \theta \|^2$ (resp. $\| b\theta\|^2$). Indeed, minimax estimation and minimax testing are essentially two different problems, see, e.g., \cite{IS_2003}, Sections 1.4.4 and 2.10. Nevertheless, ideas and methodologies in minimax estimation can inspire the construction of appropriate minimax testing procedures.\\

In the following, we focus on the construction of linear estimators based on observations from the GSM (\ref{1.0.0}). Let $\omega = (\omega_j)_{j\in \mathbb{N}}$ be a filter, i.e., a sequence taking values in the interval $[0,1]$. Then, one can estimate $\|\theta\|^2$ by the following estimator
\begin{equation}
\widehat{ \| \theta \|^2} = \sum_{j\in \mathbb{N}} \omega_j b_j^{-2} (y_j^2-\varepsilon^2).
\label{eq:linear_est1}
\end{equation}
or, in the same spirit, estimate  $\|b\theta\|^2$ by the following estimator
\begin{equation}
\widehat{\| b\theta \|^2}  = \sum_{j\in \mathbb{N}} \omega_j  (y_j^2-\varepsilon^2).
\label{eq:linear_est2}
\end{equation}
Various possible filters $\omega = (\omega_j)_{j\in \mathbb{N}}$ are available in the literature. Among them, one can mention, e.g., spectral cut-off filters (see Section \ref{subsub:sp-c-off}), Tikhonov filters, Ingster filters (see Section \ref{subsubsec:ingster-F}) or filters based on other regularization approaches. For more details regarding available regularization methods, we refer, e.g., to \cite{Bissantz_2007}, \cite{EHN_1996} and \cite{IS_2003}. \\

Having an estimator $\|\theta\|^2$ (resp. $\|b\theta\|^2$) of the form (\ref{eq:linear_est1}) (resp. (\ref{eq:linear_est2})), denoted by $\hat d$, we can construct an associated test $\Psi_{\alpha,\omega}$ as
$$ \Psi_{\alpha,\omega}= \mathbf{1}_{\lbrace \hat d > t_{\alpha,\omega} \rbrace},$$
where $t_{\alpha,\omega}$ is a threshold that (asymptotically) controls the first kind error probability $\Alpha_\varepsilon(\Psi_{\alpha,\omega})$.\\

It is important to point out at this point that, having an (asymptotic) $\alpha$-level test $\Psi_\alpha$,
\begin{itemize}
\item (non-asymptotic approach) one can try to determine the smallest possible separation radius $r_{\varepsilon,0}:=r_\varepsilon(\mathcal{E}_a,\Psi_\alpha,\beta)>0$ such that the maximal second kind error probability $\Beta_{\varepsilon}(\Theta_a(r_{\varepsilon,0}), \Psi_\alpha)$ is at most $\beta$, for any prescribed $\alpha,\beta \in ]0,1[,$
\item (asymptotic approach) one can investigate the asymptotic behavior of the maximal second kind error probability $\Beta_{\varepsilon}(\Theta_a(r_\varepsilon), \Psi_\alpha)$, for any given $r_\varepsilon >0$ and any prescribed $\alpha \in ]0,1[$.
\end{itemize}

Our aim below is, 
\begin{enumerate}
\item to construct appropriate tests that reach (at least up to a constant) the lower bounds established in Sections \ref{s:lb_nasymp} and \ref{subsubsec:ac-final-1} (Sections \ref{subsub:sp-c-off} and \ref{subsubsec:ingster-F}),
\item to bring into light hitherto unknown links between non-asymptotic and asymptotic approaches to minimax signal detection (see Section \ref{connection}).
\end{enumerate}

As mentioned previously, there exist several possible available filters. We focus below on two different kind of filters investigated in, e.g., \cite{LLM_2012} and \cite{ISS_2012}, namely, spectral cut-off and Ingster filters, respectively.

\subsubsection{Non-asymptotic control: Spectral cut-off filters}
\label{subsub:sp-c-off}
Our aim is to propose an $\alpha$-level test $\Psi_\alpha$ such that
$$ 
r_{\varepsilon}(\mathcal{E}_a,\Psi_\alpha,\beta)  \leq C r_{\varepsilon,\star},
$$ 
for some $\mathcal{C} \geq 1$, where $r_{\varepsilon,\star}>0$ has been introduced in (\ref{bi_nasymp1}).  In such a case, this will mean that lower and upper bounds for the minimax separation radius $\tilde r_{\varepsilon}>0$ match together, up to a constant.    \\

According to the previous discussion, the suggested test will be based on an estimation of $\|\theta\|^2$ (using (\ref{eq:linear_est1})). More formally, given a bandwidth $D \in \mathbb{N}$, we define
\begin{eqnarray}
\label{LLM_test}
\Psi_{D,P} &:=& \mathbf{1}_{\lbrace \sum_{j=1}^D b_j^{-2} (y_j^2-\varepsilon^2) > t_{1-\alpha,D} \rbrace} \nonumber\\
&:=& \mathbf{1}_{\lbrace T_{D,P} > t_{1-\alpha,D} \rbrace},
\end{eqnarray}
where 
\begin{equation}
\label{eq:fanis-sco-t}
T_{D,P}:= \sum_{j=1}^D b_j^{-2} (y_j^2-\varepsilon^2)
\end{equation}
and $t_{1-\alpha,D}$ denotes the $(1-\alpha)$-quantile of $T_{D,P}$ under $H_0$, i.e., the $(1-\alpha)$-quantile of the random variable $\varepsilon^2 \sum_{j=1}^D b_j^{-2} (\xi_j^2-1)$.\\

Due to the definition of $t_{1-\alpha, D}$, the spectral cut-off test $\Psi_{D,P}$ is an $\alpha$-level test. Indeed,
\begin{eqnarray*} 
\Alpha_\varepsilon(\Psi_{D,P}) := \mathbb{P}_0\left( \Psi_{D,P} =1 \right) &=& \mathbb{P}_0\left( T_{D,P} > t_{1-\alpha,D}\right) \nonumber \\&=&  \mathbb{P}\left( \varepsilon^2 \sum_{j=1}^D b_j^{-2} (\xi_j^2-1) > t_{1-\alpha,D}\right) = \alpha.
\end{eqnarray*}

Now, we turn to the control of the maximal second kind error probability $\Beta_{\varepsilon}(\Theta_a(r_\varepsilon),\Psi_{D,P})$. To this end, denote by $t_{\beta,D}(\theta)$ the $\beta$-quantile of $T_{D,P}$ under $H_1$, i.e., the term satisfying
$$ \mathbb{P}_\theta \left( T_{D,P} \leq t_{\beta,D}(\theta) \right) = \beta.$$
Then, for a given $\theta\in \mathcal{E}_a$, in order to prove that
$$ \mathbb{P}_\theta( \Psi_{D,P} =0) := \mathbb{P}_\theta\left(T_{D,P} \leq t_{1-\alpha,D} \right) = \mathbb{P}_\theta \left( \sum_{j=1}^D b_j^{-2} (y_j^2-\varepsilon^2) \leq t_{1-\alpha,D} \right) \leq \beta,$$
it suffices to show that
\begin{equation}
t_{1-\alpha,D} \leq t_{\beta,D}(\theta).
\label{eq:cond_quantile}
\end{equation}

\begin{figure}
\begin{center}
\begin{tikzpicture}[scale=1.5]
\draw [thin,->] (-0.5,0) -- (7.5,0);
\draw [pattern=vertical lines,domain=0:5.5,thick] plot (\x, {3*\x*exp(-\x)});
\path [pattern=vertical lines] (3.5,0) -- (3.5,0.1) -- (5.5,0.1) -- (5.5,0) -- (3.5,0);
\path [fill=white] (0,0) -- (3.5,0.4171) -- (3.5,0) -- (0,0);
\draw [fill=white,domain=0:3.5,thick] plot (\x, {3*\x*exp(-\x)});
\draw [domain=3:7,thick] plot (\x, {4*(\x-3)*exp(-1.2*\x+3.6)});
%\draw [domain=3:7,thick] plot (\x, {4*(\x-3)*exp(-1.2*(\x-3))});
\draw [pattern=horizontal lines,domain=3:4,thick] plot (\x, {4*(\x-3)*exp(-1.2*\x+3.6)});
\path [pattern=horizontal lines] (3,0) -- (4,1.2048) --(4,0) -- (3,0);
\draw [thick,dashed] (4,-0.5) -- (4,2.5); 
\draw [thick,dashed] (3.5, -0.5) -- (3.5,2.5);
\draw (3.1,-0.4) node {$t_{1-\alpha,D}$};
\draw (3.5,0) circle (0.03);
\draw (4.6,-0.4) node {$t_{\beta,D}(\theta)$};
\draw (4,0) circle (0.03);
\draw [thin,->] (0.9,-0.5) -- (0.9,3);
\draw (0.8,-0.2) node {$0$};
\draw [thick,->] (2.3,1.6)--(3.3,0.6);
\draw (2.2,1.7) node {\Large $\beta$};
\draw [thick,<-] (4.5,0.08) -- (5.5,1.1);
\draw (5.65,1.2) node {\Large $\alpha$};
\end{tikzpicture}
\end{center}
\caption{\textit{Illustrative comparison of the $(1-\alpha)$-quantile $t_{1_\alpha,D}$ and the $\beta$-quantile $t_{\beta,D}(\theta)$ for a fixed bandwidth $D \in \mathbb{N}$. The left hand side curve displays the density of the test statistic $T_{D,P}:=\sum_{j=1}^Db_j^{-2}(y_j^2-\epsilon^2)$, defined in (\ref{eq:fanis-sco-t}), under $H_0$, while the one on the right hand side displays the density of the same test statistics $T_{D,P}$ under $H_1$. The shaded areas are determined by the corresponding $(1-\alpha)$-quantile $t_{1-\alpha,D}$ $($vertical lines$)$ and $\beta$-quantile $t_{\beta,D}(\theta)$ $($horizontal lines$)$.}}
\label{fig:Fig3}
\end{figure}
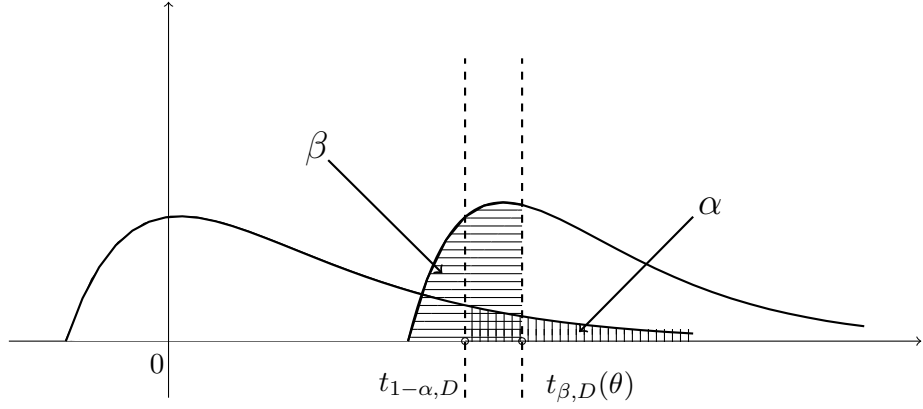

Figure \ref{fig:Fig3} provides, for a fixed bandwidth $D \in \mathbb{N}$, a heuristic illustration for the comparison (\ref{eq:cond_quantile}) between the $(1-\alpha)$-quantile $t_{1-\alpha,D}$ and the $\beta$-quantile $t_{\beta,D}(\theta)$ of the test statistic $T_{D,P}$, defined in (\ref{eq:fanis-sco-t}).  In order to compare these two terms formally, we use the following proposition.

\begin{proposition}
\label{prop:quantile}
Let $T_{D,P}$ be the test statistic defined in (\ref{eq:fanis-sco-t}), and let $t_{1-\alpha,D}$ and $t_{\beta,D}(\theta)$ be its $(1-\alpha)$-quantile under $H_0$ and $\beta$-quantile under $H_1$. Then, there exists a constant (explicitly computable) $C(\alpha)>0$ such that
$$ t_{1-\alpha,D} \leq C(\alpha) \varepsilon^2 \left( \sum_{j=1}^D b_j^{-4} \right)^{1/2},$$
and
$$  t_{\beta,D} \geq \sum_{j=1}^D \theta_j^2 - 2\sqrt{\ln(1/\beta)} \sqrt{ \varepsilon^4 \sum_{j=1}^D b_j^{-4} + 2 \varepsilon^2 \sum_{j=1}^D b_j^{-2} \theta_j^2}.$$
\end{proposition}

The proof of Proposition \ref{prop:quantile} can be found in \cite{LLM_2012} (see the construction of the upper and lower bounds in the proof of their Proposition 2). In particular, the control of $t_{1-\alpha,D}$ and $t_{\beta,D}(\theta)$ is based on deviation inequalities of appropriate independent weighted-$\chi^2$ random variables.\\

Using (\ref{eq:cond_quantile}) and Proposition \ref{prop:quantile}, one can easily see that
$$
t_{1-\alpha,D} \leq t_{\beta,D}(\theta) 
$$
if and only if
$$
\sum_{j=1}^D \theta_j^2  \geq C(\alpha) \varepsilon^2 \left( \sum_{j=1}^D b_j^{-4} \right)^{1/2} - 2\sqrt{\ln(1/\beta)} \sqrt{ \varepsilon^4 \sum_{j=1}^D b_j^{-4} + 2 \varepsilon^2 \sum_{j=1}^D b_j^{-2} \theta_j^2},
$$
which, in turn, holds true as soon as 
\begin{equation}
 \sum_{j=1}^D \theta_j^2  \geq C(\alpha,\beta) \varepsilon^2 \sqrt{\sum_{j=1}^D b_j^{-4}},
\label{eq:upper_condition}
\end{equation}
where, setting $x_\gamma=\ln(1/\gamma)$, for all $\gamma \in ]0,1[$, 
\begin{equation}
C(\alpha,\beta) =\sqrt{2x_\beta} + \sqrt{2(x_\alpha+x_\beta)} + \sqrt{2}\; \big(\sqrt{x_\alpha}+\sqrt{x_\beta}\big)^{1/2}>0
\label{clement_C}
\end{equation}
(see \cite{LLM_2012} for more details). The condition (\ref{eq:upper_condition}) ensures that
$$
\mathbb{P}_\theta(\Psi_{D,P} = 0)\leq \beta.
$$ 

The main drawback of (\ref{eq:upper_condition}) is that it is expressed in terms of a lower bound on $\sum_{j=1}^D \theta_j^2$ instead of $\|\theta\|^2$. However, since $\theta \in \mathcal{E}_a$, it follows that $\sum_{j>D} \theta_j^2 \leq a_D^{-2}$. Hence, 
\begin{equation}
\forall \;\theta \in \mathcal{E}_a, \quad  \| \theta \|^2 \geq C(\alpha,\beta) \varepsilon^2 \sqrt{\sum_{j=1}^D b_j^{-4}} + a_D^{-2}
 \quad \Rightarrow \quad \mathbb{P}_\theta(\Psi_{D,P} = 0)\leq \beta.
\label{eq:upper_condition2}
\end{equation}
Moreover, we point out that the term in the left hand side of (\ref{eq:upper_condition2}) corresponds to the sum of two antagonist quantities. Since our aim is to obtain the weaker possible bound on the energy condition, we choose a bandwidth $D:=D^{\star} \in \mathbb{N}$ such that
$$ 
\Beta_{\varepsilon}(\Theta_a(r_\varepsilon^\star),\Psi_{D^\star,P}):=
\sup_{\theta\in \Theta_a(r_{\varepsilon}^\star)}\mathbb{P}_\theta(\Psi_{D^\star,P}=0) \leq \beta, $$
where
\begin{equation}
(r_{\varepsilon}^\star)^2:= \inf_{D\in\mathbb{N}} \left[C(\alpha,\beta) \varepsilon^2 \sqrt{\sum_{j=1}^D b_j^{-4}} + a_D^{-2}\right].
\label{bi_nasymp2}
\end{equation}

Finally, thanks to (\ref{bi_nasymp1}) and (\ref{bi_nasymp2}), the values of $r_{\varepsilon,\star}$ and $r_{\varepsilon}^\star$ are coherent. In Section \ref{connection}, we show that, under some weak conditions on the sequences $a=(a_j)_{j\in \mathbb{N}}$ and $b=(b_j)_{j\in\mathbb{N}}$, non-asymptotic lower and upper bounds for the minimax separation radius $\tilde r_\varepsilon$ match together, up to a constant.

\subsubsection{Asymptotic control: Ingster  filters}
\label{subsubsec:ingster-F}
We consider a different approach, since the testing procedure will be based on an estimation of $\|b\theta\|^2$ (using (\ref{eq:linear_est2})). We will deal with a specific kind of filters which have been, to the best of our knowledge, introduced by Yuri I. Ingster in a series of papers (see, e.g., \cite{IngsterI}, \cite{IngsterII}, \cite{IngsterIII}, \cite{IS_2003}).\\

Let $\bar \theta=\bar\theta(r_\varepsilon) \in \Theta_a(r_\varepsilon)$ be the solution of the extremal problem (\ref{extrem_pb}). Then, we define the Ingster filters $\omega_{r_\varepsilon} = (\omega_{j,r_\varepsilon})_{j\in\mathbb{N}}$ as
\begin{equation} 
\omega_{j,r_\varepsilon}= \frac{b_j^2\bar \theta_j^2}{\sqrt{2\sum_{k\in \mathbb{N}} b_k^4\bar \theta_k^4}}\quad \forall \, j\in \mathbb{N}.
\label{yuri_filter}
\end{equation}
As discussed in Section \ref{s:general_test}, one can use a test of the form
\begin{eqnarray}
\Psi_{r_\varepsilon,I} 
& := & \mathbf{1}_{\lbrace \sum_{j\in\mathbb{N}} \omega_{j,r_\varepsilon} (y_j^2-\varepsilon^2) > \varepsilon^2\, t_{1-\alpha}\rbrace} \nonumber\\
& = & \mathbf{1}_{\lbrace \sum_{j\in\mathbb{N}} \omega_{j,r_\varepsilon} ((y_j/\varepsilon))^2-1) >  t_{1-\alpha}\rbrace} \nonumber \\&:=&\mathbf{1}_{\lbrace T_{r_\varepsilon,I} > t_{1-\alpha}\rbrace},
\label{yuri_test}
\end{eqnarray}
where
\begin{equation}
T_{r_\varepsilon,I} = \sum_{j\in\mathbb{N}} \omega_{j,r_\varepsilon} \left( \left(\frac{y_j}{\varepsilon}\right)^2-1 \right),
\label{test_statistics}
\end{equation}
and $t_{1-\alpha}$ denotes the $(1-\alpha)$-quantile of a standard Gaussian random variable. \\

Since $T_{r_\varepsilon,I}=\xi_\varepsilon$, where $\xi_\varepsilon$ is the quantity appeared in (\ref{fanis:loglikeasymp}), from the proof of the corresponding lower bounds, it follows that $T_{r_\varepsilon,I} \rightarrow \xi \sim {\cal N}(0,1)$ (in $\mathbb{P}_0$ probability) as $\varepsilon \rightarrow 0$. Hence, we immediately see that, 
\begin{equation}
\label{eq:alpha-asymp-yuri} 
\Alpha_\varepsilon(\Psi_{r_\varepsilon,I}) = \Phi(t_{1-\alpha}) + o_\varepsilon(1) = \alpha+o_\varepsilon(1).
\end{equation}
(This means that $\Psi_{r_\varepsilon,I}$ is, asymptotically, an $\alpha$-level test.)\\

We now consider the corresponding maximal second kind error probability $\Beta_\varepsilon(\Theta_a(r_\varepsilon),\Psi_{r_\varepsilon,I})$. The following cases are of particular interest:
\begin{itemize}
\item $u_\varepsilon(r_\varepsilon)=o_\varepsilon(1)$. According to (\ref{eq:lower_bound_zero}), $\Beta_{\varepsilon,\alpha}(\Theta_a(r_\varepsilon)) \geq 1-\alpha +o_\varepsilon(1)$. It is then impossible to distinguish $H_0$ from $H_1$, meaning that, according to Definition \ref{trivial}, we have an asymptotical trivial test. Thus, it is not needed to further study this case.
\item  $u_\varepsilon(r_\varepsilon)=O_\varepsilon(1)$. In this case, under the mild condition $\sup_{j\in\mathbb{N}} \omega_{j,r_\varepsilon}= o_\varepsilon(1)$, we establish a sharp asymptotic of Gaussian type for $\Beta_\varepsilon(\Theta_a(r_\varepsilon),\Psi_{r_\varepsilon,I})$.
\item $u_\varepsilon(r_\varepsilon) \rightarrow \infty$ as $\varepsilon \rightarrow 0$.  In this case, we establish that $\Beta_\varepsilon(\Theta_a(r_\varepsilon),\Psi_{r_\varepsilon,I}) =o_\varepsilon(1)$.
\end{itemize}

To this end, simple algebra leads to the following expressions of the expectation and the variance of the test statistics $T_{r_\varepsilon,I}$: 
\begin{equation}
\mathbb{E}_\theta[T_{r_\varepsilon,I}] = \varepsilon^{-2} \sum_{j\in \mathbb{N}} \omega_{j,r_\varepsilon} b_j^2\theta_j^2, \hspace{1cm} \mathrm{Var}_\theta[T_{r_\varepsilon,I}] = 1 + 4\varepsilon^{-2} \sum_{j\in \mathbb{N}} \omega^2_{j,r_\varepsilon} b_j^2\theta_j^2.
\label{eq:exp-var}
\end{equation}
Introduce the standardized random variable $\tilde T_{r_\varepsilon,I}$ defined as
$$ \tilde T_{r_\varepsilon,I} = \frac{T_{r_\varepsilon,I} - \mathbb{E}_\theta[T_{r_\varepsilon,I}]}{\sqrt{\mathrm{Var}_\theta[T_{r_\varepsilon,I}]}},$$
where the $\mathbb{E}_\theta[T_{r_\varepsilon,I}]$ and $\mathrm{Var}_\theta[T_{r_\varepsilon,I}]$ have been computed in (\ref{eq:exp-var}). Define 
\begin{equation} 
\omega_{0,r_\varepsilon}:=\sup_{j\in\mathbb{N}}\frac{b_j^2\bar \theta_j^2}{\sqrt{2\sum_{k\in \mathbb{N}} b_k^4\bar \theta_k^4}} :=\sup_{j\in\mathbb{N}} \omega_{j,r_\varepsilon},
\label{eq:omega_0}
\end{equation}
where $\bar \theta=\bar\theta(r_\varepsilon) \in \Theta_a(r_\varepsilon)$ is the extremal sequence, i.e., the solution of the extremal problem (\ref{extrem_pb}). (Note that, using (\ref{eq:omega_0}), $1 \leq \mathrm{Var}_\theta[T_{r_\varepsilon,I}] \leq 1 + 4\,\omega_{0,r_\varepsilon}\, \mathbb{E}_\theta[T_{r_\varepsilon,I}]$.)\\

In order to proceed, we need the following proposition.

\begin{proposition}
\label{p:function_h}
Let $T_{r_\varepsilon,I}$ be the test statistic defined in (\ref{test_statistics})
Let $h(r_\varepsilon,\theta):=\mathbb{E}_\theta[T_{r_\varepsilon,I}]$, where $\mathbb{E}_\theta[T_{r_\varepsilon,I}]$ is computed in (\ref{eq:exp-var}). Then
$$\inf_{\theta \in \Theta_a(r_\varepsilon)} h(r_\varepsilon,\theta)  = u_\varepsilon(r_\varepsilon).$$
\end{proposition}

The proof of Proposition  \ref{p:function_h} can be found in \cite{ISS_2012}, supplementary material, Lemma 11.1. \\

\noindent
{\bf Case 1} ($u_\varepsilon(r_\varepsilon)=O_\varepsilon(1)$) Following the proof of Theorem 4.1 of \cite{ISS_2012}, supplementary material, Section 11.1, using Lyapunov's conditions and (\ref{eq:exp-var}), it follows that, as soon as $\omega_{0,r_\varepsilon} =o_\varepsilon(1)$, 
\begin{itemize}
\item $\tilde T_{r_\varepsilon,I}$ is asymptotically standard Gaussian under $\mathbb{P}_\theta$.
\item $\mathrm{Var}_\theta[T_{r_\varepsilon,I}] =1+o_\varepsilon(1)$.
\end{itemize}
Therefore, we get that 
$$ \mathbb{P}_\theta(\Psi_{r_\varepsilon,I}=0) = \mathbb{P}_\theta \left(\tilde T_{r_\varepsilon,I} \leq \frac{t_{1-\alpha} - \mathbb{E}_\theta[T_{r_\varepsilon,I}]}{\sqrt{\mathrm{Var}_\theta[T_{r_\varepsilon,I}]}} \right)= \Phi(t_{1-\alpha} - \mathbb{E}_\theta[T_{r_\varepsilon,I}])+o_\varepsilon(1).$$
Using Proposition \ref{p:function_h}, we arrive at
\begin{eqnarray}
\Beta_\varepsilon(\Theta_a(r_\varepsilon),\Psi_{r_\varepsilon,I}) 
& = & \sup_{\theta \in \Theta_a(r_\varepsilon)}  \mathbb{P}_\theta(\Psi_{r_\varepsilon,I}=0) \nonumber \\
& = & \sup_{\theta \in \Theta_a(r_\varepsilon)}  \Phi(t_{1-\alpha} - \mathbb{E}_\theta[T_{r_\varepsilon,I}])+o_\varepsilon(1)  \nonumber\\
& := & \sup_{\theta \in \Theta_a(r_\varepsilon)}  \Phi(t_{1-\alpha} - h(r_\varepsilon,\theta))+o_\varepsilon(1) \nonumber\\
& = &   \Phi\left(t_{1-\alpha} - \inf_{\theta \in \Theta_a(r_\varepsilon)} h(r_\varepsilon,\theta) \right)+o_\varepsilon(1) \nonumber \\
& = &  \Phi\left(t_{1-\alpha} - u_\varepsilon(r_\varepsilon) \right)+o_\varepsilon(1).
\label{sharp_asymp1}
\end{eqnarray}
Therefore, using (\ref{sharp_asymp1}), we get that
\begin{eqnarray}
\Beta_{\varepsilon,\alpha}(\Theta_a(r_\varepsilon)) &\leq& \Beta_\varepsilon(\Theta_a(r_\varepsilon),\Psi_{r_\varepsilon,I})  \nonumber \\ &=& \Phi(t_{1-\alpha} - u_\varepsilon(r_\varepsilon)) +o_\varepsilon(1).
\label{sharp_asymp2}
\end{eqnarray}
\noindent
(Note that, in the particular case that $r_\varepsilon >0$ satisfies
$u_\varepsilon(r_\varepsilon)=t_{1-\alpha} - t_\beta$, it is immediately seen that   
$\Beta_{\varepsilon,\alpha}(\Theta_a(r_\varepsilon)) \leq \beta + o_\varepsilon(1)$.)

\vspace{0.2cm}
\noindent
{\bf Case 2} ($u_\varepsilon(r_\varepsilon) \rightarrow +\infty$ as $\varepsilon \rightarrow 0$) 
Using Proposition \ref{p:function_h} and (\ref{eq:exp-var}), it follows that, for all $\theta \in \Theta_a(r_\varepsilon)$,
$$\mathbb{E}_\theta[T_{r_\varepsilon,I}] :=h(r_\varepsilon,\theta)  \geq  \inf_{\theta \in \Theta_a(r_\varepsilon)} h(r_\varepsilon,\theta)  = u_\varepsilon(r_\varepsilon) \rightarrow + \infty \quad \mbox{as} \quad 
\varepsilon \rightarrow 0.
$$
Therefore, using Markov's inequality,
\begin{eqnarray}
\label{eq:Markov-ineq}
\Beta_\varepsilon(\Theta_a(r_\varepsilon),\Psi_{r_\varepsilon,I}) &:=& \sup_{\theta\in \Theta_a(r_\varepsilon)} \mathbb{P}_\theta(\Psi=0) \nonumber\\
&=& \sup_{\theta\in \Theta_a(r_\varepsilon)} \mathbb{P}_\theta \left(\tilde T_{r_\varepsilon,I} \leq \frac{t_{1-\alpha} - \mathbb{E}_\theta[T_{r_\varepsilon,I}]}{\sqrt{\mathrm{Var}_\theta[T_{r_\varepsilon,I}]}} \right) \nonumber\\
&\leq& \sup_{\theta\in \Theta_a(r_\varepsilon)} \mathbb{P}_\theta \left(|\tilde T_{r_\varepsilon,I}| \geq \frac{\mathbb{E}_\theta[T_{r_\varepsilon,I}]-t_{1-\alpha}}{\sqrt{\mathrm{Var}_\theta[T_{r_\varepsilon,I}]}} \right) \nonumber\\
&\leq& \sup_{\theta\in \Theta_a(r_\varepsilon)} \frac{\mathrm{Var}_\theta[T_{r_\varepsilon,I}]}{(t_{1-\alpha} - \mathbb{E}_\theta[T_{r_\varepsilon,I}])^2} \nonumber\\
&\leq& \sup_{\theta\in \Theta_a(r_\varepsilon)} \frac{1 + 4\,\omega_{0,r_\varepsilon}\, h(r_\varepsilon,\theta)}{(h(r_\varepsilon,\theta) -t_{1-\alpha})^2}\nonumber\\
&\sim& \frac{1}{\inf_{\theta \in \Theta_a(r_\varepsilon)} h(r_\varepsilon,\theta)}:=\frac{1}{u_\varepsilon(r_\varepsilon)}=o_\varepsilon(1).
\end{eqnarray}

\begin{remark}
\label{rem:123}
{\rm It is worth mentioning that one cannot determine at this point  the radius $r_{\varepsilon,2}>0$ (considered in the general strategy of Section \ref{asympt_general} for constructing upper bounds). This more or less amounts to solve the equation $u_\varepsilon(r_{\varepsilon,2}) = \mathcal{O}_\varepsilon(1)$ as $\varepsilon \rightarrow 0$. 
This cannot be accomplished unless the sequences $a=(a_j)_{j\in \mathbb{N}}$ and $b=(b_j)_{j\in\mathbb{N}}$ are explicitly given. We refer again to, e.g., \cite{ISS_2012} for more details or to Section \ref{s:example} where an example of a mildly ill-posed inverse problem is treated for illustrative purposes.}
\end{remark}

Below, we first formalize the results for the lower and upper bounds presented above and explain their meaning for practical purposes (Section \ref{subsec:gen_res}). We then bring into light hitherto unknown links between non-asymptotic and asymptotic approaches to minimax signal detection (Section \ref{subsec:explicit_seq} and Section \ref{s:example}).

\section{Connections between non-asymptotic and asymptotic frameworks}
\label{connection}

\subsection{General Results}
\label{subsec:gen_res}
In Section \ref{s:control}, lower and upper bounds on the minimax separation radius $\tilde r_{\varepsilon}$ and the minimax second kind error probability $\Beta_{\varepsilon,\alpha}(\Theta_a(r_\varepsilon))$ were independently treated. In the following theorems, these results are gathered in unified manners.\\

We first focus our attention to the non-asymptotic paradigm.

\begin{thm}
\label{thm:1}
(Non-asymptotic framework)  Assume that $Y=(Y_j)_{j\in\mathbb{N}}$ are observations from the GSM (\ref{1.0.0}), and consider the signal detection problem (\ref{testing_pb1}) with ${\cal F}$ defined in (\ref{def:f-set}). Let $\alpha, \beta \in\, ]0,1[$ be given. Then, for every $\varepsilon >0$, the minimax separation radius $\tilde r_{\varepsilon}$ is controlled by 
\begin{equation}
\sup_{D\in \mathbb{N}} \left[ c(\alpha,\beta) \varepsilon^2 \sqrt{\sum_{j=1}^D b_j^{-4}} \wedge a_D^{-2} \right] \leq  \tilde r^2_{\varepsilon} \leq \inf_{D\in\mathbb{N}} \left[C(\alpha,\beta) \varepsilon^2 \sqrt{\sum_{j=1}^D b_j^{-4}} + a_D^{-2}\right], 
\label{eq:control_tilde_r}
\end{equation} 
where the constants $c(\alpha,\beta)$ and $C(\alpha,\beta)$ are respectively given in (\ref{clement_c}) and (\ref{clement_C}).
\end{thm}

In order to shed some light on the meaning of (\ref{eq:control_tilde_r}), the following comments are in order:

\begin{itemize}
\item One cannot ensure that both lower and upper bounds on the minimax separation radius $\tilde r_\varepsilon$ in (\ref{eq:control_tilde_r}) match together, unless (weak) conditions on the sequences $a=(a_j)_{j\in\mathbb{N}}$ and $b=(b_j)_{j\in\mathbb{N}}$ are at hand. A discussion on that point is provided below (see Theorem \ref{thm:2}). 
\item Nevertheless, we point out that these bounds are coherent since they involve the same quantities, namely, $a_D^{-2}$ and $\varepsilon^2 \sqrt{\sum_{j=1}^D b_j^{-4}}$, for any given bandwidth $D\in \mathbb{N}$,  as well as positive constants $c(\alpha,\beta)$ and $C(\alpha,\beta)$, depending on $\alpha$ and $\beta$ only.  
\item A careful look into the discussion presented in the previous section indicates that the term $a_D^{-2}$ can be related to a (in fact an upper bound on the) `bias' term in the sense that it measures the amount of signal that is missed using the spectral cut-off test $\Psi_{D,P}$ (see, e.g., (\ref{eq:upper_condition})). Recall that  the sequence $a=(a_j)_{j\in\mathbb{N}}$ characterizes the smoothness of the underlying signal $\theta$. Obviously, the smoother the signal of interest, the easier the testing problem in the sense that the minimax separation radius $\tilde r_\varepsilon$ becomes smaller. 
\item In the same spirit, $\varepsilon^2 \sqrt{\sum_{j=1}^D b_j^{-4}}$ can be related to a `standard deviation' term that corresponds to the estimation of the term $\|\theta\|^2$ using the spectral cut-off test $\Psi_{D,P}$. When $b_j=1$, for all $j\in\mathbb{N}$, (i.e., the direct problem) this term is of order $\varepsilon^2 \sqrt{D}$. This particular case has been discussed in detail in \cite{Baraud}, Section 3. On the other hand, the case when $b_j \rightarrow 0 $ as $j\rightarrow +\infty$, corresponds to ill-posed inverse problems. In this case, the signal detection problem becomes harder in the sense that the minimax separation radius $\tilde r_\varepsilon$ strongly depends on the decay of the sequence $b=(b_j)_{j\in\mathbb{N}}$ towards $0$ and becomes larger than the corresponding one in the direct problem. 
\end{itemize}

In summary, in order to precisely compute the minimax separation radius $\tilde r_\varepsilon$, explicit sequences of $a=(a_j)_{j\in\mathbb{N}}$ and $b=(b_j)_{j\in\mathbb{N}}$ are needed to control the trade-off between the two antagonist terms, i.e., the `bias' and the `standard deviation' terms, $a_D^{-2}$ and $\varepsilon^2 \sqrt{\sum_{j=1}^D b_j^{-4}}$, respectively. This will be elaborated on Section \ref{s:example} below, where an example of a mildly ill-posed inverse problem is used for illustrative purposes.\\

We now turn our attention to the asymptotic paradigm.

\begin{thm}
\label{thm:1a} (Asymptotic framework) Assume that $Y=(Y_j)_{j\in\mathbb{N}}$ are observations from the GSM (\ref{1.0.0}), and consider the signal detection problem (\ref{testing_pb1}) with ${\cal F}$ defined in (\ref{def:f-set}). Let a radius $r_\varepsilon>0$ be fixed, and let $\alpha \in\, ]0,1[$ be given. Let $u_\varepsilon(r_\varepsilon)$ and \,$\omega_{0,r_\varepsilon}$ denote the solution of the extremal problem defined in (\ref{extrem_sol}) and the term introduced in (\ref{eq:omega_0}), respectively.   
\begin{enumerate}
\item[\rm{(a)}] If $$u_\varepsilon(r_\varepsilon)=o_\varepsilon(1),$$
then 
$$
\Beta_{\varepsilon,\alpha}(\Theta_a(r_\varepsilon)) = 1-\alpha +o_\varepsilon(1).
$$  
\item[\rm{(b)}] If  $$u_\varepsilon(r_\varepsilon)=\mathcal{O}_\varepsilon(1) \quad \text{and} \quad \omega_{0,r_\varepsilon}=o_\varepsilon(1),
$$
then
$$ \Beta_{\varepsilon,\alpha}(\Theta_a(r_\varepsilon)) = \Phi(t_{1-\alpha}-u_\varepsilon(r_\varepsilon)) + o_\varepsilon(1).
$$  
\item[\rm{(c)}] If  $$u_\varepsilon(r_\varepsilon) \rightarrow +\infty \quad  \text{as} \quad \varepsilon \rightarrow 0,
$$
then
$$ \Beta_{\varepsilon,\alpha}(\Theta_a(r_\varepsilon)) =o_\varepsilon(1).$$  
\end{enumerate}
\end{thm}

It is evident from Theorem \ref{thm:1a} that the minimax signal detection problem in the asymptotic framework essentially reduces to the study of the extremal problem (\ref{extrem_pb}). Indeed, the corresponding solution given in (\ref{extrem_sol}) governs both the lower and the upper bounds on the minimax second kind error probability $\Beta_{\varepsilon,\alpha}(\Theta_a(r_\varepsilon))$. The three different regimes mentioned in Theorem \ref{thm:1a} are of particular interest and require at this step some additional explanations: 

\begin{itemize}
\item If $u_\varepsilon(r_\varepsilon)=o_\varepsilon(1)$, then, according to Definition \ref{trivial}, an asymptotical non-trivial minimax hypothesis testing problem is not possible. In other words, it is impossible to distinguish between $H_0$ and $H_1$.
\item If $u_\varepsilon(r_\varepsilon)=\mathcal{O}_\varepsilon(1)$ and $\omega_{0,r_\varepsilon}=o_\varepsilon(1)$,  then one can precisely describe the shape of the minimax second kind error probability $\Beta_{\varepsilon,\alpha}(\Theta_a(r_\varepsilon))$ since it possesses a sharp asymptotic of Gaussian type. It is also evident that $\Beta_{\varepsilon,\alpha}(\Theta_a(r_\varepsilon))\in \,]0,1-\alpha[$. This means that the minimax signal detection problem is asymptotically non-trivial (i.e., $\Beta_{\varepsilon,\alpha}(\Theta_a(r_\varepsilon))>0$) but that $H_0$ and $H_1$ can be asymptotically always separated (i.e., $\Beta_{\varepsilon,\alpha}(\Theta_a(r_\varepsilon)) <1-\alpha$). Note that in this particular case that $u_\varepsilon(r_\varepsilon)=\mathcal{O}_\varepsilon(1)$ and $\omega_{0,r_\varepsilon}=o_\varepsilon(1)$, the Ingster test $\Psi_{r_\varepsilon,I}$ defined in (\ref{yuri_test}) is asymptotical minimax according to Definition \ref{def:am}.
\item If $u_\varepsilon(r_\varepsilon)\rightarrow +\infty$ as $\varepsilon \rightarrow 0$, then the minimax second kind error probability $\Beta_{\varepsilon,\alpha}(\Theta_a(r_\varepsilon))=o_\varepsilon(1)$. In particular, the test $\Psi_{r_\varepsilon,I}$ constructed in (\ref{yuri_test})-(\ref{test_statistics}) asymptotically always separates $H_0$ from $H_1$.   
\end{itemize} 

\begin{remark}
{\rm Theorem \ref{thm:1a} does not treat the case where
\begin{equation}
\label{eq:non_treated}
u_\varepsilon(r_\varepsilon)=\mathcal{O}_\varepsilon(1) \quad \text{and} \quad \omega_{0,r_\varepsilon}\not\rightarrow 0  \quad  \text{as} \quad \varepsilon \rightarrow 0.
\end{equation}
In such a case, the lower bound (\ref{sharp_asymp}) is still valid but can be, in fact, improved by showing that
$$
\mathrm{liminf}_{\varepsilon \rightarrow 0}\, \Beta_{\varepsilon,\alpha}(\Theta_a(r_\varepsilon)) >1-\alpha \quad \mbox{for any} \quad \alpha \in ]0,1[,
$$
i.e., the minimax signal detection problem is asymptotically trivial (see the proof of Theorem 4.1 of \cite{ISS_2012}, supplementary material, Section 11.1.). It is worth pointing out at this point that if (\ref{eq:non_treated}) holds, then the minimax second kind error probability $\Beta_{\varepsilon,\alpha}(\Theta_a(r_\varepsilon))$ asymptotically belongs to the set $\{0,1-\alpha\}$, for any $\alpha \in ]0,1[$, depending on the behavior of any given $r_\varepsilon >0$.   

The case where  $\omega_{0,r_\varepsilon}\not\rightarrow 0$  as $\varepsilon \rightarrow 0$ exists in, e.g., the case of severely ill-posed inverse problems with the class of analytic functions (super-smooth functions), i.e., $b_j \asymp e^{-jt}$, $j \in \mathbb{N}$, for some $t>0$, and $a_j \asymp e^{js}$, $j \in \mathbb{N}$, for some $s>0$, respectively. Indeed,
$$
\omega_{0,r_\varepsilon}:=\sup_{j \in \mathbb{N}}\frac{b_j^2\bar
\theta_j^2}{\sqrt{2\sum_{k \in \mathbb{N}}b_k^4\bar \theta_k^4}}\sim \frac{z_0^2
e^{-2 t m}}{z_0^2 e^{-2 t m}}\asymp 1 \not\to 0, \quad  \text{as} \quad \varepsilon \rightarrow 0,
$$
for some quantities $z_0 \in \mathbb{R}$ and $m \in [1, \infty)$  (see
Theorem 4.3 and Remark 4.4 in \cite{ISS_2012}). We also refer to Section \ref{subsec:ingster_nonasy} below for a similar computation in a mildly ill-posed inverse problem setting.}
\end{remark}

\begin{remark}
\label{rem:111}
{\rm Theorem \ref{thm:1a} does not provide an immediate expression for the minimax separation rate $\bar r_\varepsilon$. In practice, however, both terms $r_{\varepsilon,1}$ and $r_{\varepsilon,2}$ required in the construction of the lower and upper bounds, respectively, sketched in Section \ref{asympt_general}, are derived from the same equation:  $u_\varepsilon(r_{\varepsilon,1})=u_\varepsilon(r_{\varepsilon,2})=\mathcal{O}_\varepsilon(1)$. Then one can, `in general', check the implications
$$
\frac{r_\varepsilon}{r_{\varepsilon,1}} \rightarrow 0 \; \Rightarrow \; u_\varepsilon(r_\varepsilon)=o_\varepsilon(1) \quad \mbox{and} \quad  \frac{r_\varepsilon}{r_{\varepsilon,2}} \rightarrow +\infty \; \Rightarrow \; u_\varepsilon(r_\varepsilon) \rightarrow +\infty, \quad  \text{as} \quad \varepsilon \rightarrow 0,
$$ 
which, thanks to Theorem \ref{thm:1a} and Definition \ref{def:separation}, allows one to conclude.
As mentioned previously, this task cannot be accomplished unless explicit expressions on the sequences $a=(a_j)_{j\in\mathbb{N}}$ and $b=(b_j)_{j\in\mathbb{N}}$ are given. Explicit calculation of the minimax separation rate $\bar r_\varepsilon$ in a mildly ill-posed inverse problem is provided in Section \ref{s:asymp_rate}.}
\end{remark}

The proofs of the assertions in Theorem \ref{thm:1} and Theorem \ref{thm:1a} are direct consequences of the discussion provided in Section \ref{s:control}, concerning the control of the upper and lower bounds, for both minimax separation radius $\tilde r_{\varepsilon}$ and maximal second kind error probability $\Beta_{\varepsilon,\alpha}(\Theta_a(r_\varepsilon))$. Detailed arguments  and related discussions can be found in, e.g., \cite{Baraud}, \cite{ISS_2012} and \cite{LLM_2012}.

\subsection{Deriving the minimax separation rate $\bar r_\varepsilon$ from bounds on the minimax separation radius $\tilde r_{\varepsilon}$.}
\label{subsec:explicit_seq}

The following theorem shows that, under some mild conditions on the growth of the sequences $a=(a_j)_{j\in \mathbb{N}}$ and $b^{-1}=(b^{-1}_j)_{j\in\mathbb{N}}$, one can derive the asymptotic order of the minimax separation rate $\bar r_\varepsilon$ from the bounds on the minimax separation radius $\tilde r_{\varepsilon}$ given in (\ref{eq:control_tilde_r}). 

\begin{proposition}
\label{thm:2}
Assume that $Y=(Y_j)_{j\in\mathbb{N}}$ are observations from the GSM (\ref{1.0.0}), and consider the signal detection problem (\ref{testing_pb1}) with ${\cal F}$ defined in (\ref{def:f-set}).  
Assume that both sequences $a=(a_j)_{j\in \mathbb{N}}$ and $b^{-1}=(b_j^{-1})_{j\in\mathbb{N}}$ are non-decreasing and that they satisfy
\begin{equation}
a_\star \leq \frac{a_{D-1}}{a_{D}} \leq a^\star \quad \mathrm{and} \quad b_\star \leq \frac{b_{D-1}}{b_{D}} \leq b^\star \quad \mbox{for all} \quad D > 1,
\label{eq:hyp_ab}
\end{equation}
for some constants $0<a_\star \leq a^\star <\infty$ and $0<b_\star \leq b^\star <\infty$.
Let $\alpha, \beta \in\, ]0,1[$ be given. Then, there exists a constant $C\geq 1$ such that
$$ \inf_{D\in\mathbb{N}} \left[C(\alpha,\beta) \varepsilon^2 \sqrt{\sum_{j=1}^D b_j^{-4}} + a_D^{-2}\right] \leq C \sup_{D\in \mathbb{N}} \left[ c(\alpha,\beta) \varepsilon^2 \sqrt{\sum_{j=1}^D b_j^{-4}} \wedge a_D^{-2} \right] ,$$
where the constants $c(\alpha,\beta)$ and $C(\alpha,\beta)$ are respectively given in (\ref{clement_c}) and (\ref{clement_C}). In particular, both lower and upper bounds in (\ref{eq:control_tilde_r}) are of the same order.

%the minimax separation radius $\tilde r_{\varepsilon}$ and the minimax separation rate $\bar r_\varepsilon$ are, up to a constant, asymptotically equivalent, i.e.,
%\begin{equation}
%\label{Fanis-Equiv-Rates-0}  
%\frac{\tilde r_{\varepsilon}}{\bar r_\varepsilon} = \mathcal{O}_\varepsilon(1). 
%\end{equation}
\end{proposition}

The control (\ref{eq:control_tilde_r}) hence proposes a sharp description of the minimax separation radius $\tilde r_\epsilon$ as soon as (\ref{eq:hyp_ab}) is satisfied. Since $\bar r_\epsilon / \tilde r_\epsilon = O_\epsilon(1)$, the asymptotic minimax separation rate $\bar r_\epsilon$ can thus be determine from (\ref{eq:control_tilde_r}). On the other hand, a non-asymptotic bound that matches asymptotic known results can be considered as (rate) optimal. Hence, although the motivations differ, both asymptotic and non-asymptotic approaches provide a similar description on the minimax signal detection problem at hand.

\begin{remark}
{\rm \bigskip
We note also that the condition (\ref{eq:hyp_ab}) is satisfied for various combinations of interest, among them: (i) mildly ill-posed inverse problems ($b_j \asymp j^{-t}$, $j\in\mathbb{N}$, for some $t>0$) with ordinary smooth functions ($a_j  \asymp j^{s}$, $j\in\mathbb{N}$, for some $s>0$), (ii) severely ill-posed inverse problems ($b_j \asymp e^{-jt}$, $j\in\mathbb{N}$, for some $t>0$) with ordinary smooth functions ($a_j \asymp j^{s}$, $j\in\mathbb{N}$, for some $s>0$), and (iii) mildly ill-posed inverse problems ($b_j \asymp j^{-t}$, $j\in\mathbb{N}$, for some $t>0$) with super-smooth functions ($a_j \asymp e^{js}$, $j\in\mathbb{N}$, for some $s>0$). Among the possible situations where the condition (\ref{eq:hyp_ab}) is not satisfied, one can mention, for instance, power-exponential behaviors ($a_j \asymp e^{j^l s}$, $j\in\mathbb{N}$, for some $s>0$ and $l>1$, or $b_j \asymp e^{-j^rt}$, $j\in\mathbb{N}$, for some $t>0$ and $r>1$).}
\label{rem:2-F}
\end{remark}

\noindent
\textsc{Proof}. 
%First, we prove that the lower and the upper bounds in (\ref{eq:control_tilde_r}) are of the same order. To this end, l
Let the bandwidth $D_0 \in \mathbb{N}$ satisfy
\begin{equation}
D_0 = \mathrm{arg} \sup_{D\in \mathbb{N}} \left[ c(\alpha,\beta) \varepsilon^2 \sqrt{\sum_{j=1}^D b_j^{-4}} \wedge a_D^{-2} \right].
\label{eq:D_0}
\end{equation}
We restrict ourselves to the following case
$$ a_{D_0}^{-2} \leq c(\alpha,\beta) \varepsilon^2 \sqrt{\sum_{j=1}^{D_0} b_j^{-4}}.$$
(The other case follows similarly along the same lines of proof.) Then, thanks to (\ref{eq:hyp_ab}), we get 
\begin{eqnarray*}
\inf_{D\in\mathbb{N}} \left[C(\alpha,\beta) \varepsilon^2 \sqrt{\sum_{j=1}^D b_j^{-4}} + a_D^{-2}\right]
& \leq & C(\alpha,\beta) \varepsilon^2 \sqrt{\sum_{j=1}^{D_0} b_j^{-4}} + a_{D_0}^{-2}\\
& \leq & (C(\alpha,\beta) + c(\alpha,\beta))  \varepsilon^2 \sqrt{\sum_{j=1}^{D_0} b_j^{-4}}\\
& \leq & C \varepsilon^2 \sqrt{\sum_{j=1}^{D_0-1} b_j^{-4}},
\end{eqnarray*}
for some constant $C>0$ that can be explicitly computed. Note that
$$ a_{D_0}^{-2} \leq c(\alpha,\beta) \varepsilon^2 \sqrt{\sum_{j=1}^{D_0} b_j^{-4}} 
\quad
\text{implies}
\quad
a_{D_0-1}^2 > c(\alpha,\beta) \varepsilon^2 \sqrt{\sum_{j=1}^{D_0-1} b_j^{-4}},$$
since, otherwise, we arrive at a contradiction, due to the definition of $D_0 \in \mathbb{N}$ in (\ref{eq:D_0}). Hence,
\begin{eqnarray*}
\inf_{D\in\mathbb{N}} \left[C(\alpha,\beta) \varepsilon^2 \sqrt{\sum_{j=1}^D b_j^{-4}} + a_D^{-2}\right]
& \leq & C \varepsilon^2 \sqrt{\sum_{j=1}^{D_0-1} b_j^{-4}},\\
& \leq & C \left[ \varepsilon^2 \sqrt{\sum_{j=1}^{D_0-1} b_j^{-4}} \wedge a_{D_0-1}^{-2} \right],\\
& \leq & C  \sup_{D\in \mathbb{N}} \left[ c(\alpha,\beta) \varepsilon^2 \sqrt{\sum_{j=1}^D b_j^{-4}} \wedge a_D^{-2} \right].
\end{eqnarray*}
In other words, there exists some constant $C\geq 1$ such that
\begin{eqnarray}
\lefteqn{\sup_{D\in \mathbb{N}} \left[ c(\alpha,\beta) \varepsilon^2 \sqrt{\sum_{j=1}^D b_j^{-4}} \wedge a_D^{-2} \right]} \nonumber \\
& \leq & \inf_{D\in\mathbb{N}} \left[C(\alpha,\beta) \varepsilon^2 \sqrt{\sum_{j=1}^D b_j^{-4}} + a_D^{-2}\right] \leq  C  \sup_{D\in \mathbb{N}} \left[ c(\alpha,\beta) \varepsilon^2 \sqrt{\sum_{j=1}^D b_j^{-4}} \wedge a_D^{-2} \right]. \label{eq:order}
\end{eqnarray}
Hence, the lower and the upper bounds in (\ref{eq:control_tilde_r}) are of the same order. This concludes the proof of the proposition. 

\begin{flushright}
$\Box$
\end{flushright}

\begin{remark}
{\rm  According to Definition \ref{def:am}, and as soon as (\ref{eq:hyp_ab}) is satisfied, the spectral cut-off test $\Psi_{D,P}$ defined in (\ref{LLM_test}), with bandwidth $D:=D_0 \in \mathbb{N}$ selected as in (\ref{eq:D_0}), is asymptotical minimax consistent. Indeed, given $\theta \in \mathcal{E}_a$ and a radius $r_\varepsilon$ such that $\| \theta \| \geq r_\varepsilon$,
\begin{eqnarray*}
\mathbb{P}_{\theta}(\Psi_{D_0,P} =0) 
& = & \mathbb{P}_{\theta}\left( T_{D_0,P} \leq t_{1-\alpha,D_0} \right), \\ 
& = & \mathbb{P}_{\theta}\left( T_{D_0,P}- \mathbb{E}[T_{D_0,P}] \leq t_{1-\alpha,D_0}  - \sum_{j=1}^{D_0} \theta_j^2 \right),  \\ 
& \leq & \mathbb{P}_{\theta}\left( \left| \mathbb{E}[T_{D_0,P}] - T_{D_0,P} \right|  \geq \sum_{j=1}^{D_0} \theta_j^2 - t_{1-\alpha,D_0}  \right), \\ 
& \leq &  \frac{\varepsilon^4 \sum_{j=1}^{D_0} b_j^{-4}}{ \left( \sum_{j=1}^{D_0} \theta_j^2 - t_{1-\alpha,D_0}  \right)^2 }.
\end{eqnarray*}
Then, since $\bar r_\varepsilon/\tilde r_\varepsilon =O_\varepsilon(1)$ as $\varepsilon \rightarrow 0$, using Proposition \ref{prop:quantile}, we get 
\begin{eqnarray*}
\sum_{j=1}^{D_0} \theta_j^2 - t_{1-\alpha,D_0} 
& \geq & \| \theta \|^2 - C\left( \varepsilon^2 \sqrt{ \sum_{j=1}^{D_0} b_j^{-4}} + a_{D_0}^{-2} \right),\\
& \geq &  r^2_\varepsilon - C \tilde r^2_\varepsilon,\\
& \geq &  r^2_{\varepsilon}(1-o_\varepsilon(1)),
\end{eqnarray*}	
as soon as $r_\varepsilon/\bar r_\varepsilon \rightarrow +\infty$ as $\varepsilon \rightarrow 0$.  Finally, we obtain 
$$ \mathbb{P}_{\theta}(\Psi_{D_0,P} =0)  \leq \frac{C \bar r_\varepsilon^4}{r_\varepsilon^4(1-o_\varepsilon(1))} =o_\varepsilon(1),$$
which entails
$$ \beta_\varepsilon(\Theta_a(r_\varepsilon),\Psi_{D_0,P}) = o_\varepsilon(1) \quad \mathrm{if} \quad \frac{r_\varepsilon}{\bar r_\varepsilon} \rightarrow +\infty \quad \mathrm{as} \quad \varepsilon \rightarrow 0.$$
}
\end{remark}

\subsection{An illustrative example: a mildly ill-posed inverse problem}
\label{s:example}

Our aim below is to illustrate the results presented in Section \ref{subsec:gen_res} and Section \ref{subsec:explicit_seq}. To this end, we address the minimax signal detection problem of a mildly ill-posed inverse problem. Namely, we will assume that 
\begin{equation}
a_j \asymp j^{s}, \;\; \text{for some}\;\; s>0, \quad \mathrm{and} \quad b_j \asymp j^{-t}, \;\; \text{for some}\;\; t>0, \;\; \text{for all}\;\;  j\in \mathbb{N}.
\label{eq:mipp}
\end{equation}
Our aim in this context is multifold:
\begin{itemize}
\item First, we consider an asymptotic analysis of the minimax separation radius $\tilde r_\varepsilon$ based on the inequality (\ref{eq:control_tilde_r}).
\item Second, we explicitly compute the minimax separation rate $\bar r_\varepsilon$ though a careful analysis of the extremal problem (\ref{extrem_sol}).
\item Third, we provide a non-asymptotic  analysis of the Ingster test $\Psi_{r_\varepsilon,I}$, defined in (\ref{yuri_test}). In particular, we compute its associated separation radius $r_\varepsilon(\mathcal{E}_a,\Psi_{r_\varepsilon,I},\beta)$ and show that, up to constant, it coincides with the minimax separation radius $\tilde r_\varepsilon$. 
\item Fourth, we present an asymptotic analysis of the spectral cut-off test $\Psi_{D,P}$, defined in (\ref{LLM_test}). In particular, for an appropriate bandwidth $D:=\tilde{D} \in \mathbb{N}$, we prove that the maximal second kind error probability $\Beta_{\varepsilon}(\Theta_a(r_\varepsilon),\Psi_{\tilde D,P})$,  is asymptotically bounded from above by a quantity that possesses a Gaussian shape.
\end{itemize}

These results demonstrate that tools used to derive `asymptotic' results can be exploited to draw `non-asymptotic' conclusions, and vice-versa.

\subsubsection{Asymptotic analysis of the minimax separation radius $\tilde r_{\varepsilon}$} 

We are interested in the asymptotic behavior of the minimax separation radius $\tilde r_{\varepsilon}$. Recall from Theorem \ref{thm:1} that, for any $\varepsilon >0$, 
$$ \sup_{D\in \mathbb{N}} \left[ c(\alpha,\beta) \varepsilon^2 \sqrt{\sum_{j=1}^D b_j^{-4}} \wedge a_D^{-2} \right] \leq  \tilde r^2_{\varepsilon} \leq \inf_{D\in\mathbb{N}} \left[C(\alpha,\beta) \varepsilon^2 \sqrt{\sum_{j=1}^D b_j^{-4}} + a_D^{-2}\right].$$ 
 Moreover, according to (\ref{eq:order}), both the upper and the lower bounds in the above inequality are of the same order. Indeed, the constraint (\ref{eq:hyp_ab}) of Theorem \ref{thm:2} is satisfied in the setting (\ref{eq:mipp}).  Hence, we are now able to characterize the asymptotic value of the minimax separation radius $\tilde r_{\varepsilon}$. \\

Simple algebra shows that  
$$ \sum_{j=1}^D b_j^{-4} = \sum_{j=1}^D j^{4t} = C D^{4t+1}(1+o(1)) \quad \mathrm{as} \quad D\rightarrow +\infty,$$
for some constant $C>0$. Hence, the bandwidth $D_0 \in \mathbb{N}$, introduced in (\ref{eq:D_0}), satisfies 
$$ \varepsilon^2 \sqrt{\sum_{j=1}^{D_0} b_j^{-4} } =\mathcal{O}( a_{D_0}^{-2})  \quad \text{iff} \quad\varepsilon^2 D_0^{2t+1/2}  = \mathcal{O}(D_0^{-2s}) \quad \text{iff} \quad D_0 = \mathcal{O}_\varepsilon\left( \varepsilon^{\frac{-2}{2s+2t+1/2}} \right).$$
We then deduce from the previous computation that the minimax separation radius $\tilde r_{\varepsilon}$ satisfies
\begin{equation}
\label{Fanis-Equiv-Rates-1}  \tilde r^2_{\varepsilon} = \mathcal{O}(D_0^{-2s}) = \mathcal{O}_\varepsilon \left( \varepsilon^{\frac{4s}{2s+2t+1/2}} \right).
\end{equation}

\subsubsection{Computation of the minimax separation rate $\bar r_\varepsilon$}
\label{s:asymp_rate}
Following Remark \ref{rem:111}, an explicit computation of the function $r_\varepsilon \mapsto u_\varepsilon(r_\varepsilon)$ is required in order to retrieve the minimax separation rate $\bar r_\varepsilon$ from the solution of the equation $u_\varepsilon(r_\varepsilon)= \mathcal{O}_\varepsilon (1)$.\\

We first need to solve the extremal problem (\ref{extrem_sol}) defined as
\begin{equation}\label{eq0-FF}
u^2_\varepsilon(r_\varepsilon) = \frac{1}{2\varepsilon^4} \inf_{\theta \in \Theta_a(r_\varepsilon)} \sum_{j\in \mathbb{N}} b_j^{4}\theta_j^4. 
\end{equation}
This problem is solved via Lagrange multipliers. In particular, the extremal sequence, i.e., the solution of the above mentioned extremal problem, appears to be of the form 
$$ \bar \theta^2_j(r_\varepsilon) = z_0^2 b_j^{-4} (1-\mathcal{A}a_j^2)_+, \quad j \in \mathbb{N},$$ 
where the quantities $z_0:=z_{0,\varepsilon}$ and $\mathcal{A}:=\mathcal{A}_{\varepsilon}$ are determined by
the equations
\begin{equation}\label{eq1-FF}
\sum_{j\in\mathbb{N}} b_j^{-2} \bar \theta^2_j(r_\varepsilon) \quad \mbox{and} \quad
\sum_{j\in\mathbb{N}} a_j^2 b_j^{-2} \bar \theta^2_j(r_\varepsilon)=1.
\end{equation}

\begin{remark}
{\rm The quantity $\mathcal{A}$ determines the so-called {\em efficient dimension}
$m$ in specific ill-posed inverse problems: since
$a_j$ is an increasing sequence, the
efficient dimension is the quantity $m=m_\varepsilon\in [1,\infty)$ such that
$\mathcal{A}a_{[m]}^2\leq 1<\mathcal{A}a_{[m]+1}^2$, see, e.g., \cite{ISS_2012}, supplementary material, Section 11. Moreover, a unique solution to the system of equations (\ref{eq1-FF}) exists for $r_\varepsilon >0$ small enough, due to the fact that $\sum_{j \in \mathbb{N}}b_j^{-1}=+\infty$ (see,  Proposition 11.2 of \cite{ISS_2012}, supplementary material, Section 11).}
\end{remark}

The equations (\ref{eq0-FF})-(\ref{eq1-FF}) are immediately rewritten in the form
\begin{equation} 
\left\lbrace
\begin{array}{l}
r_\varepsilon^2 = z_0^2 J_1,\\
1 = z_0^2 \mathcal{A}^{-1} J_2,
\end{array}
\right.
\label{eq:extremal1}
\end{equation}
with
\begin{eqnarray*}
J_1 & = & \sum_{j\in \mathbb{N}} b_j^{-4} (1-\mathcal{A}a_j^2)_+,\\
J_2 & = & \mathcal{A}\sum_{j\in \mathbb{N}} a_j^2 b_j^{-4} (1-\mathcal{A}a_j^2)_+.
\end{eqnarray*}
In particular,  the extremal problem (\ref{extrem_sol}) takes the form
\begin{equation}
u_\varepsilon^2(r_\varepsilon) = \varepsilon^{-4} z_0^4 J_0/2, \quad \mathrm{where} \quad J_0 = J_1-J_2 = \sum_{j\in \mathbb{N}} b_j^{-4} (1-\mathcal{A}a_j^2)_+^2.
\label{eq:extremal4}
\end{equation}
Setting $\mathcal{R}=\mathcal{A}^{-1/2s}$, simple computations lead to 
\begin{eqnarray*}
J_1 & = &  \sum_{j\in \mathbb{N}} b_j^{-4} (1-\mathcal{A}a_j^2)_+,\\
& = &  \sum_{j:\, a_j^2\leq \mathcal{A}^{-1}} b_j^{-4} -  \mathcal{A} \sum_{j:\, a_j^2\leq \mathcal{A}^{-1}} b_j^{-4}a_j^2,\\
& = &  \sum_{j:\, j\leq \mathcal{R}} j^{4t} -  \mathcal{A} \sum_{j:\, j\leq \mathcal{R}} j^{4t+2s},\\
& = & \mathcal{C}_1  \mathcal{R}^{4t+1}(1+o(1)) \quad \mathrm{as} \quad \mathcal{R}\rightarrow + \infty.
\end{eqnarray*}
Using similar algebra, one can prove that
$$ J_2=\mathcal{C}_2  \mathcal{R}^{4t+1}(1+o(1)) \quad \mathrm{and} \quad J_0=\mathcal{C}_0  \mathcal{R}^{4t+1}(1+o(1)) \quad \mathrm{as} \quad \mathcal{R}\rightarrow + \infty.$$ 
In particular, we get from (\ref{eq:extremal1}) that
\begin{equation}
r_\varepsilon^2 = \mathcal{A} \frac{J_1}{J_2} = \mathcal{R}^{-2s} \frac{\mathcal{C}_1}{\mathcal{C}_2}(1+o(1))\;\; \text{as} \;\; \mathcal{R} \rightarrow +\infty.
\label{eq:comput1}
\end{equation}
Therefore, combining the above results,
\begin{eqnarray}
u_\varepsilon(r_\varepsilon) &= &\left( \frac{r_\varepsilon}{\varepsilon} \right)^4 \frac{J_0}{2J_1^2} 
 \nonumber \\
 &=& \left( \frac{r_\varepsilon}{\varepsilon} \right)^4 \frac{\mathcal{C}_0  \mathcal{R}^{4t+1}}{2 \mathcal{C}^2_2  \mathcal{R}^{2(4t+1)}}(1+o(1)) \quad \text{as} \;\; \mathcal{R} \rightarrow +\infty, \nonumber\\
&=&  \left( \frac{r_\varepsilon}{\varepsilon} \right)^4  \frac{r^{-(4t+1)/s}_\varepsilon}{r^{-(8t+2)/s}_\varepsilon} \mathcal{O}(1) \quad  \text{as} \;\; r_\varepsilon \rightarrow 0, \nonumber\\
&=& \mathcal{O} \left( \varepsilon^{-4}  r^{(4s+4t+1)/s}_\varepsilon \right)
\label{eq:comput_0} \quad  \text{as} \;\; r_\varepsilon \rightarrow 0.
\end{eqnarray}

The expression in (\ref{eq:comput_0}) provides an explicit form for the function $r_\varepsilon \mapsto u_\varepsilon(r_\varepsilon)$ that is required in order to retrieve the minimax separation rate $\bar r_\varepsilon$ from the solution of the equation $u_\varepsilon(r_\varepsilon)= \mathcal{O}_\varepsilon (1)$. 
Our next task is to solve this equation. Using (\ref{eq:comput_0}), we immediately get
\begin{equation}
u_\varepsilon(r_\varepsilon) = \mathcal{O}_\varepsilon (1) \quad \text{iff} \quad  \varepsilon^{-4}  r^{(4s+4t+1)/s}_\varepsilon =\mathcal{O}_\varepsilon (1) \quad \text{iff} \quad r_\varepsilon  = \mathcal{O}_\varepsilon \left( \varepsilon^{\frac{2s}{2s+2t+1/2}} \right).
\label{eq:comput222}
\end{equation} 

In order to conclude our discussion, we need to prove that the minimax separation rate $\bar r_\varepsilon$ is of the following order
\begin{equation}
\label{eq:fff}
r_{\varepsilon,0} = \mathcal{O}_\varepsilon \left( \varepsilon^{\frac{2s}{2s+2t+1/2}} \right).
\end{equation}
To this end, we remark that, for any $r_\varepsilon >0$,
\begin{itemize}
\item If $r_\varepsilon/ r_{\varepsilon,0} \rightarrow 0$ then, using (\ref{eq:comput_0}), it is easily seen that $u_\varepsilon(r_\varepsilon) =o_\varepsilon(1)$. Hence, according to Theorem \ref{thm:1a}, $\Beta_\varepsilon ( \Theta_a(r_\varepsilon)) = 1-\alpha +o_\varepsilon(1)$.
\item If $r_\varepsilon/ r_{\varepsilon,0} \rightarrow +\infty$ then, using (\ref{eq:comput_0}), it is easily seen that $u_\varepsilon(r_\varepsilon) \rightarrow +\infty$, as $\varepsilon \rightarrow 0$. Hence, according to Theorem \ref{thm:1a}, $\Beta_\varepsilon ( \Theta_a(r_\varepsilon)) = o_\varepsilon(1)$.
\end{itemize}
Therefore, Definition \ref{def:separation} allows to conclude that $r_{\varepsilon,0}$ in (\ref{eq:fff}) is indeed the minimax separation rate $\bar r_\varepsilon$. (Note that, in view of (\ref{Fanis-Equiv-Rates-1}) and (\ref{eq:fff}), the minimax separation radius $\tilde r_\varepsilon$ and the minimax separation rate $\bar r_\varepsilon$ are of the same asymptotic order, as expected according to previous discussion.)

\subsubsection{Non-asymptotic analysis of the Ingster test $\Psi_{r_\varepsilon,I}$} 
\label{subsec:ingster_nonasy}
We present a non-asymptotic study of the Ingster test $\Psi_{r_\varepsilon,I}$, defined in (\ref{yuri_test}). We show that the statistical performances of the Ingster test $\Psi_{r_\varepsilon,I}$ and the spectral cut-off test $\Psi_{D,P}$ defined in (\ref{LLM_test}) are comparable. In particular, the Ingster test $\Psi_{r_\varepsilon,I}$ appears to be powerful in the sense of Definition \ref{def:powerFanis}, namely, 
$$\ r_\varepsilon(\mathcal{E}_a,\Psi_{r_\varepsilon,I},\beta) \leq \mathcal{C} \tilde r_{\varepsilon},$$
for a fixed $\beta \in ]0,1[$ and some constant $\mathcal{C} \geq 1$, for an appropriately selected radius $r_{\varepsilon} >0$. 

\begin{proposition}
\label{prop:NA}
Let $\alpha, \beta \in ]0,1[$ be given. Define 
\begin{equation}
\rho_\varepsilon^2 :=\inf_{{\cal R} \geq 2} [\mathcal{C}_0 \mathcal{R}^{-2s} \vee c'(\alpha,\beta) \varepsilon^2 \mathcal{R}^{2t+1/2}],
\label{eq:rho}
\end{equation}
for some positive constants $\mathcal{C}_0 $ and $c'(\alpha,\beta)$ than can be explicitly computed. Let $ \Psi_I^\star := \Psi_{\rho_\varepsilon,I},$
where $\Psi_{.,I}$ is the Ingster test defined in (\ref{yuri_test}), and let $\rho_\varepsilon>0$ be the radius defined in (\ref{eq:rho}). Then, there exists constants $\mathcal{C} \geq 1$ and $\varepsilon_0>0$ such that, for all $0<\varepsilon < \varepsilon_0$, the separation radius of $\Psi_I^\star$
satisfies
\begin{equation} 
r_{\varepsilon}(\mathcal{E}_a, \Psi_{I}^\star,\beta) \leq \mathcal{C} \tilde r_{\varepsilon}.
\label{eq:yuri-nonn-asymp-11}
\end{equation}
\end{proposition}

\noindent 
\textsc{Proof of Proposition \ref{prop:NA}}. Let $r_\varepsilon >0$ be a given radius. Using the same arguments as in (\ref{eq:Markov-ineq}), we get
$$
\Beta_\varepsilon(\Theta_a(r_\varepsilon), \Psi_{r_\varepsilon,I}) \leq \sup_{\theta\in \Theta_a(r_\varepsilon)} \frac{1 + 4\,\omega_{0,r_\varepsilon}\, h(r_\varepsilon,\theta)}{(h(r_\varepsilon,\theta) -t_{1-\alpha})^2}.
$$
Then, there exists an explicit constant $C_{\alpha,\beta} >0$ such that
\begin{eqnarray}
u_\varepsilon(r_\varepsilon) \geq C_{\alpha,\beta} \; & \Rightarrow & h(r_\varepsilon,\theta) \geq C_{\alpha,\beta}, \quad \forall \; \theta\in \Theta_a(r_\varepsilon) \nonumber\\
& \Rightarrow & \frac{1 + 4\,\omega_{0,r_\varepsilon}\, h(r_\varepsilon,\theta)}{(h(r_\varepsilon,\theta) -t_{1-\alpha})^2} \leq  \beta \quad \forall \; \theta \in \Theta_a(r_\varepsilon) \nonumber\\
& \Rightarrow & \sup_{\theta\in \Theta_a(r_\varepsilon)} \frac{1 + 4\,\omega_{0,r_\varepsilon}\, h(r_\varepsilon,\theta)}{(h(r_\varepsilon,\theta) -t_{1-\alpha})^2} \leq \beta \nonumber\\
& \Rightarrow & \Beta_\varepsilon(\Theta_a(r_\varepsilon), \Psi_{r_\varepsilon,I}) \leq \beta.
\end{eqnarray}

Our task now is to find a condition on $r_\varepsilon>0$ that will guarantee the validity of the above inequality $u_\varepsilon(r_\varepsilon) \geq C_{\alpha,\beta}$. Working along the lines of Section \ref{s:asymp_rate}, we then arrive at
$$ r_\varepsilon^2 = z_0^2 J_1, \quad \mathcal{A} = z_0^2 J_2 \quad \mathrm{and} \quad u_\varepsilon^2(r_\varepsilon) = \varepsilon^{-4}r_\varepsilon^4 \frac{J_0}{2J_1^2}.$$
Hence, we see that
\begin{equation} 
u_\varepsilon(r_\varepsilon) \geq C_{\alpha,\beta} \quad \text{iff} \quad r_\varepsilon^2 \geq \sqrt{2}\, C_{\alpha,\beta}\, \varepsilon^2\, \frac{J_1}{J_0^{1/2}}.
\label{eq:inter1}
\end{equation}
Moreover, 
$$ J_1 =  \sum_{j\in \mathbb{N}} b_j^{-4} (1-\mathcal{A}a_j^2)_+ \leq \sum_{j:j\leq \mathcal{R}} j^{4t} \leq C_1' \mathcal{R}^{4t+1},$$
and
$$ J_1  \geq \sum_{j:j\leq \mathcal{R}} j^{4t} - \mathcal{R}^{-2s} \sum_{j:j\leq \mathcal{R}} j^{4t + 2s} \geq C_1 \mathcal{R}^{4t+1},$$
for all $\mathcal{R}\geq 2$, where $\mathcal{R}=\mathcal{A}^{-1/2s}$, and for some positive constants $C_1, C_1'$ (depending on $s$ and $t$ only). In the same spirit, we can also prove that, for all $\mathcal{R}\geq 2$,
$$ C_2 \mathcal{R}^{4t+1} \leq J_2 \leq  C_2' \mathcal{R}^{4t+1} \quad \mathrm{and} \quad C_0 \mathcal{R}^{4t+1} \leq J_0 \leq  C_0' \mathcal{R}^{4t+1},$$
for some positive constants $C_0,C_0',C_1,C_1'$ (depending on $s$ and $t$ only). Hence, we get that
\begin{equation} 
\sqrt{2} C_{\alpha,\beta} \varepsilon^2 \frac{J_1}{J_0^{1/2}} \leq c'(\alpha,\beta) \varepsilon^2 \mathcal{R}^{2t+1/2},
\label{eq:inter2}
\end{equation}
for some constant $c'(\alpha,\beta) >0$. Therefore, we deduce from (\ref{eq:inter1})-(\ref{eq:inter2}), that
$$ 
r_\varepsilon^2 \geq c'(\alpha,\beta) \varepsilon^2 \mathcal{R}^{2t+1/2}
\quad \Rightarrow \quad u_\varepsilon(r_\varepsilon) \geq C_{\alpha,\beta}.$$
Using the same kind of algebra, we get from (\ref{eq:comput1}) that $r_\varepsilon^2 = \mathcal{R}^{-2s} J_1/J_0$. Hence, for all ${\cal R} \geq 2,$ 
$$ \mathcal{C}_0 \mathcal{R}^{-2s} \leq r_\varepsilon^2 \leq \mathcal{C}_1 \mathcal{R}^{-2s},$$
for some positive constants $\mathcal{C}_0, \mathcal{C}_1$.
Finally, for all $\mathcal{R} \geq 2$,
$$ 
r_\varepsilon^2 \geq  \mathcal{C}_0 \mathcal{R}^{-2s} \vee c'(\alpha,\beta) \varepsilon^2 \mathcal{R}^{2t+1/2}
\quad \Rightarrow \quad
u_\varepsilon(r_\varepsilon) \geq C_{\alpha,\beta} \quad \Rightarrow  \quad \Beta_\varepsilon(\Theta_a(r_\varepsilon), \Psi_{r_\varepsilon,I}) \leq \beta.$$
Hence, taking $r_\varepsilon := \rho_\varepsilon$, where $\rho_\varepsilon$ is defined in (\ref{eq:rho}), we immediately get that
$$  \Beta_\varepsilon(\Theta_a(\rho_\varepsilon), \Psi_I^\star) \leq \beta,$$
which implies 
$$ r_{\varepsilon}(\mathcal{E}_a, \Psi_{I}^\star,\beta) \leq \rho_\varepsilon.$$

To conclude, it suffices to show that there exists a constant $\mathcal{C} \geq 1$ and $\varepsilon_0>0$ such that, for all $0<\varepsilon < \varepsilon_0$,
$ \rho_\varepsilon \leq \mathcal{C} \tilde r_{\varepsilon}.$ This, however, holds true working along the lines of the proof of (\ref{eq:order}). This concludes the proof of the proposition. \begin{flushright} $\Box$ \end{flushright}

Concerning Proposition \ref{prop:NA}, the following comments are in order:
\begin{itemize}
\item
The considered Ingster test $\Psi_I^\star$, designed for asymptotic purposes, can be, somehow, employed in the non-asymptotic framework. It appears, that we recover existing non-asymptotic upper bounds, namely, for all $0 < \varepsilon <\varepsilon_0$, the Ingster test $\Psi_I^\star$ is powerful according to Definition \ref{def:powerFanis}. The value $\varepsilon_0>0$ guarantees that the optimal bandwidth $D_0 \in \mathbb{N}$ in (\ref{eq:D_0}) satisfies the requirement $2 \leq D_0 <+\infty$ which, in turn, ensures that $\rho_\varepsilon$ and $\tilde r_\varepsilon$ are, indeed, of the same order. 
\item
The term $\rho_\varepsilon$ involved in the construction of $\Psi_I^\star:=\Psi_{I,\rho_\varepsilon}$ plays the role of a tuning (regularization) parameter. In a sense, the parameter $\rho_\varepsilon$ plays a similar role to the bandwidth $D^\star$ in (\ref{bi_nasymp2}). Hence, it provides a trade-off between the two competing terms `bias' and `standard deviation' involved in (\ref{eq:rho}).     
\item
As we have seen in (\ref{eq:alpha-asymp-yuri}), the Ingster test $\Psi_{r_\varepsilon,I}$ is an asymptotic $\alpha$-level test for all $r_\varepsilon >0$. Hence, a non-asymptotic control of the first kind error probability $\Alpha_\varepsilon(\Psi_I^\star)$ would be necessary in order to provide a fully non-asymptotic treatment for the Ingster test $\Psi_I^\star$. This can be easily accomplished by replacing the $(1-\alpha)$-quantile $t_{1-\alpha}$ of a standard Gaussian random variable in (\ref{yuri_test}) by an appropriate $(1-\alpha)$-quantile. Then, the upper bound (\ref{eq:yuri-nonn-asymp-11}) presented in Proposition \ref{prop:NA} still holds true, up to some constants. 
\end{itemize}

\subsubsection{Asymptotics of Gaussian type for the spectral cut-off test $\Psi_{D,P}$} 
\label{subsec:spectral_nonasy}
To conclude, we present an asymptotic analysis of the spectral cut-off test $\Psi_{D,P}$, defined in (\ref{LLM_test}). In particular, as $\epsilon \rightarrow 0$, we prove that the maximal second kind error probability 
$\Beta_{\varepsilon}(\Theta_a(r_\varepsilon),\Psi_{\tilde D,P})$, for an appropriate bandwidth $\tilde D \in \mathbb{N}$, is asymptotically bounded from above by a quantity that possesses a Gaussian shape.

\begin{proposition}
\label{prop:spectral_nonasy}
Let $\alpha \in ]0,1[$ be given. Let $r_\varepsilon >0$ be a radius satisfying $u_\varepsilon(r_\varepsilon)=\mathcal{O}_\varepsilon(1)$. Let also $\Psi_{\tilde D,P}$ be the spectral cut-off test defined in (\ref{LLM_test}) with bandwidth $\tilde D  \in \mathbb{N}$ satisfying
\begin{equation}
\label{eq:Dstar_value}
\tilde D := \arg\max_{D \in \mathbb{N}} \left\{C(\alpha) \varepsilon^2 \sqrt{\sum_{j=1}^{D} b_j^{-4}} +a_{D}^{-2} \leq \frac{r^2_\varepsilon}{2}\right\},
\end{equation}
where the positive constant $C(\alpha)$ is defined in Proposition \ref{prop:quantile}.
Then, for any sequence $h_\varepsilon \in ]0,1[$ satisfying $h_\varepsilon=o_\varepsilon(1)$ and $h_\varepsilon \tilde D^{1/4} \rightarrow +\infty$ as $\varepsilon \rightarrow 0$,
\begin{equation}
\label{eq:gaussian_shape_cuttoff-0}
\Beta_{\varepsilon}(\Theta_a(r_\varepsilon),\Psi_{\tilde D,P}) \leq \Phi\left(t_{1-\alpha}-(1-h_\varepsilon) \frac{a_{\tilde D}^{-2}}{\varepsilon^2 \sqrt{\sum_{j=1}^{\tilde D} b_j^{-4}}}\right) + o_\varepsilon(1).
\end{equation}
\end{proposition}

\noindent 
\textsc{Proof of Proposition \ref{prop:spectral_nonasy}} Consider the spectral cut-off test $\Psi_{\tilde D,P}$ defined in (\ref{LLM_test}), with bandwidth $\tilde D \in \mathbb{N}$ selected as in (\ref{eq:Dstar_value}). For any $\theta\in \Theta_a(r_\varepsilon)$ and any sequence $h_\varepsilon \in ]0,1[$ (that will be made precise later on)
\begin{eqnarray}
\mathbb{P}_\theta(\Psi_{\tilde D,P}=0) &=& \mathbb{P}_\theta\left(T_{\tilde D,P} \leq 
t_{1-\alpha,\tilde D}\right) \nonumber\\
&=&\mathbb{P}_\theta\left( \sum_{j=1}^{\tilde D} b_j^{-2} (y_j^2-\varepsilon^2) \leq t_{1-\alpha,\tilde D}\right) \nonumber\\
&=& \mathbb{P}_\theta\left(\varepsilon^2 \sum_{j=1}^{\tilde D} b_j^{-2} (\xi_j^2-1) + 2\varepsilon \sum_{j=1}^{\tilde D} b_j^{-1}\theta_j \xi_j \leq t_{1-\alpha,\tilde D} -
\sum_{j=1}^{\tilde D}\theta_j^2\right) \nonumber\\
&\leq &  \mathbb{P}_\theta\left(\varepsilon^2 \sum_{j=1}^{\tilde D} b_j^{-2} (\xi_j^2-1) \leq (1-h_\varepsilon)\bigg(t_{1-\alpha,\tilde D} -
\sum_{j=1}^{\tilde D}\theta_j^2\bigg)\right) \nonumber\\
&& \;\; + \,
\mathbb{P}_\theta\left(2\varepsilon \sum_{j=1}^{\tilde D} b_j^{-1}\theta_j \xi_j \leq h_\varepsilon \bigg(t_{1-\alpha,\tilde D} -
\sum_{j=1}^{\tilde D}\theta_j^2\bigg)\right) \nonumber\\
&:=& T_1 + T_2,
\label{eq:T1-T2}
\end{eqnarray}
where, for the last inequality, we used the fact that, for any $t \in \mathbb{R}$ and any random variables $X$ and $Y$,
$$
\{ X+Y \leq t\} \subseteq \{ X \leq (1-h_\varepsilon)t\} \cup \{Y \leq h_\varepsilon t\}.
$$

Below, our aim is 
\begin{itemize}
\item to show that, asymptotically, $T_1$ has a Gaussian shape of the form (\ref{eq:gaussian_shape_cuttoff}),
\item to prove that $T_2=o_\varepsilon(1)$,
\item to study the asymptotic behavior of the threshold $t_{1-\alpha,\tilde D}$.
\end{itemize}
\medskip

\noindent
{\bf Control of $T_1$:} For any $\delta >0$, simple algebra shows that, for any bandwidth $D \in \mathbb{N}$,
$$
\frac{\varepsilon^{2(2+\delta)} \sum_{j=1}^{D} b_j^{-2(2+\delta)} \mathbb{E}|\xi_j^2-1|^{2+\delta}}{\left(\varepsilon^2 \sqrt{\sum_{j=1}^{D} b_j^{-4}}\right)^{2+\delta}} \asymp \frac{D^{2(2+\delta)t+1}}{D^{(4t+1)(2+\delta)/2}} \asymp D^{-\delta/2}=o(1) \quad \mbox{as} \quad D \rightarrow +\infty.
$$
Hence, by Lyapunov's condition, 
$$
\frac{\varepsilon^2 \sum_{j=1}^{\tilde D} b_j^{-2} (\xi_j^2-1)}{\varepsilon^2 \sqrt{\sum_{j=1}^{\tilde D} b_j^{-4}}} \stackrel{\mathcal{L}}{\longrightarrow} \mathcal{N}(0,1) \quad \mbox{as} \quad \varepsilon \rightarrow 0 \quad (\mbox{since}\quad \tilde D \rightarrow +\infty).
$$
Then, it follows that
\begin{equation}
\label{eq:T1-control-1}
T_1= \Phi \left( \frac{(1-h_\varepsilon)\bigg(t_{1-\alpha,\tilde D} -
\sum_{j=1}^{\tilde D}\theta_j^2\bigg)}{\varepsilon^2 \sqrt{\sum_{j=1}^{\tilde D} b_j^{-4}}} \right) + o_\varepsilon(1).
\end{equation}
\vspace{0.05cm}

\noindent
{\bf Control of $T_2$:}  Since $(\xi{_j})_{j \in \mathbb{N}}$ are independent standard Gaussian random variables,
\begin{eqnarray}
T_2&:=&\mathbb{P}_\theta\left(2\varepsilon \sum_{j=1}^{\tilde D} b_j^{-1}\theta_j \xi_j \leq h_\varepsilon \bigg(t_{1-\alpha,\tilde D} -
\sum_{j=1}^{\tilde D}\theta_j^2\bigg)\right) \nonumber\\
&=& \mathbb{P}_{\theta} \left(Z \leq \frac{h_\varepsilon}{2} \frac{\bigg(t_{1-\alpha,\tilde D} -
\sum_{j=1}^{\tilde D}\theta_j^2\bigg)}{\varepsilon \sqrt{\sum_{j=1}^{\tilde D} b_j^{-2} \theta_j^2}} \right) \quad (\mbox{where}\;\;Z \sim \mathcal{N}(0,1)).
\label{eq:T2-control-110}
\end{eqnarray}
Then, according to (\ref{eq:Dstar_value}), for any $\theta\in \Theta_a(r_\varepsilon)$,
\begin{eqnarray}
\frac{\sum_{j=1}^{\tilde D}\theta_j^2-t_{1-\alpha,\tilde D}}{\varepsilon \sqrt{\sum_{j=1}^{\tilde D} b_j^{-2} \theta_j^2}} &\geq& \frac{\|\theta\|^2-t_{1-\alpha,\tilde D}-a_{\tilde D}^{-2}}{\varepsilon\, (\max_{1 \leq j \leq \tilde D} b_j^{-1}) \|\theta\|} \geq
\frac{1}{C}\left( \frac{\|\theta\|^2-r^2_{\varepsilon}/2}{r_{\varepsilon}\|\theta\|}\right) \tilde D^{1/4}\nonumber\\
&\geq&\frac{\tilde D^{1/4}}{2C} \frac{\|\theta\|}{r_{\varepsilon}} \nonumber\\
&\geq& \frac{1}{2C}\tilde D^{1/4} \rightarrow +\infty \quad \mbox{as} \quad \varepsilon \rightarrow 0 \quad (\mbox{since}\quad \tilde D \rightarrow +\infty),
\label{eq:T2-control-22}
\end{eqnarray}
where for the second inequality we used the fact that
$$
\varepsilon \, \bigg(\max_{1 \leq j \leq \tilde D} b_j^{-1}\bigg) \asymp \varepsilon \tilde D^{t} \asymp \frac{\left(\varepsilon^2 \sqrt{\sum_{j=1}^{\tilde D} b_j^{-4}}\right)^{1/2}}{\tilde D^{1/4}} \leq C \frac{r_\varepsilon}{\tilde D^{1/4}},
$$
for some constant $C>0$. Hence, using (\ref{eq:T2-control-110}) and (\ref{eq:T2-control-22}), it follows that, as soon as $h_\varepsilon \tilde D^{1/4} \rightarrow +\infty$ as $\varepsilon \rightarrow 0$,
\begin{equation}
\label{eq:T2-control-1}
T_2=o_\varepsilon(1).
\end{equation}
\vspace{0.05cm}

\noindent
{\bf Behavior of $t_{1-\alpha,\tilde D}$:}  First we show that
\begin{equation}
\label{contrl_quant}
(1-h_\varepsilon)\frac{t_{1-\alpha,\tilde D}}{\varepsilon^2 \sqrt{\sum_{j=1}^{\tilde D} b_j^{-4}}}  = t_{1-\alpha} +  o_\varepsilon(1).
\end{equation}
Indeed, according to the definition of $t_{1-\alpha,\tilde D}$,
\begin{eqnarray}
\label{eq:t1-d-size}
&&\mathbb{P}_0\left( \sum_{j=1}^{\tilde D} b_j^{-2} (y_j^2-\varepsilon^2) \leq t_{1-\alpha,\tilde D}\right) = 1-\alpha \nonumber\\
&\Leftrightarrow& \mathbb{P}\left( \frac{\varepsilon^2 \sum_{j=1}^{\tilde D} b_j^{-2} (\xi_j^2-1)}{\varepsilon^2 \sqrt{\sum_{j=1}^{\tilde D} b_j^{-4}}} \leq \frac{t_{1-\alpha,\tilde D}}{\varepsilon^2 \sqrt{\sum_{j=1}^{\tilde D} b_j^{-4}}}\right) =1-\alpha \nonumber\\
&\Leftrightarrow& \Phi^{-1}_{\varepsilon}(1-\alpha)=\frac{t_{1-\alpha,\tilde D}}{\varepsilon^2 \sqrt{\sum_{j=1}^{\tilde D} b_j^{-4}}}, \nonumber
\end{eqnarray}
where, for any $s \in \mathbb{R}$,
$$
\Phi_{\varepsilon}(s):=\mathbb{P}\left( \frac{\varepsilon^2 \sum_{j=1}^{\tilde D} b_j^{-2} (\xi_j^2-1)}{\varepsilon^2 \sqrt{\sum_{j=1}^{\tilde D} b_j^{-4}}} \leq s \right).
$$
Then, as above, using the Central Limit Theorem with Lyapunov's condition and Lemma 21.2 in \cite{vdv_1998}, we get
\begin{eqnarray}
&&\Phi_\varepsilon(s) \rightarrow \Phi(s) \quad \mbox{as} \quad \varepsilon \rightarrow 0 \quad (\forall \;s \in \mathbb{R}) \nonumber\\
&\Leftrightarrow& \Phi^{-1}_\varepsilon(u) \rightarrow \Phi^{-1}(u) \quad \mbox{as} \quad \varepsilon \rightarrow 0 \quad (\forall \;u \in ]0,1[)
\end{eqnarray}
In particular, for any $\alpha \in ]0,1[$,
\begin{eqnarray}
& & \Phi^{-1}_\varepsilon(1-\alpha) \rightarrow \Phi^{-1}(1-\alpha) \quad \mbox{as} \quad \varepsilon \rightarrow 0 \nonumber\\
&\Leftrightarrow& \frac{t_{1-\alpha,\tilde D}}{\varepsilon^2 \sqrt{\sum_{j=1}^{\tilde D} b_j^{-4}}}  = t_{1-\alpha} +  o_\varepsilon(1).
\end{eqnarray}
Finally, taking into account that $h_\varepsilon=o_\varepsilon(1)$, (\ref{contrl_quant}) holds true.\\

\noindent
{\bf Completing the proof:} Using (\ref{eq:T1-control-1}) and (\ref{contrl_quant}), it follows that
\begin{equation}
\label{eq:T2-control-11}
T_1= \Phi \left( t_{1-\alpha} - \frac{(1-h_\varepsilon)
\sum_{j=1}^{\tilde D}\theta_j^2}{\varepsilon^2 \sqrt{\sum_{j=1}^{\tilde D} b_j^{-4}}} + o_\varepsilon(1) \right) + o_\varepsilon(1).
\end{equation}
According to (\ref{eq:Dstar_value}), for any $\theta\in \Theta_a(r_\varepsilon)$,
\begin{eqnarray}
(1-h_\varepsilon) \frac{\sum_{j=1}^{\tilde D}\theta_j^2}{\varepsilon^2 \sqrt{\sum_{j=1}^{\tilde D} b_j^{-4}}}
&=& (1-h_\varepsilon) \frac{\|\theta\|^2-\sum_{j=1}^{\tilde D} \theta_j^{2}}{\varepsilon^2 \sqrt{\sum_{j=1}^{\tilde D} b_j^{-4}}} \nonumber\\
&\geq& (1-h_\varepsilon) \frac{r_\varepsilon^2-a_{\tilde D}^{-2}}{\varepsilon^2 \sqrt{\sum_{j=1}^{\tilde D} b_j^{-4}}} \nonumber\\
&\geq& (1-h_\varepsilon) \frac{a_{\tilde D}^{-2}}{\varepsilon^2 \sqrt{\sum_{j=1}^{\tilde D} b_j^{-4}}}. 
\end{eqnarray}
Hence, using (\ref{eq:T1-T2}), (\ref{eq:T1-control-1}), (\ref{eq:T2-control-1}) and the Mean Value Theorem,
\begin{eqnarray}
\Beta_{\varepsilon}(\Theta_a(r_\varepsilon),\Psi_{\tilde D,P}) &=&
\sup_{\theta\in \Theta_a(r_\varepsilon)} \mathbb{P}_\theta(\Psi_{\tilde D,P}=0) \nonumber\\
&\leq&
 \Phi\left(t_{1-\alpha}-(1-h_\varepsilon) \frac{a_{\tilde D}^{-2}}{\varepsilon^2 \sqrt{\sum_{j=1}^{\tilde D} b_j^{-4}}}\right) + o_\varepsilon(1).\nonumber
\end{eqnarray}
Hence, (\ref{eq:gaussian_shape_cuttoff-0}) holds true, and this completes the proof of the proposition.  \begin{flushright} $\Box$ \end{flushright}

\begin{figure}
\begin{center}
\begin{tikzpicture}[scale=1.60]
\draw [thick,->] (-0.5,0) -- (5,0);
\draw (5.25,0) node {$r_\varepsilon$};
\draw [thick, ->] (0,-0.5) -- (0,4);
\draw (-0.15,-0.15) node {\large $0$};
\draw (-0.2,4.3) node {{\large $\beta_{\varepsilon,\alpha}(\Theta_a(r_\varepsilon))$}};
\draw [thick] (-0,3.3) -- (4.5,3.3);
\draw (-0.25,3.3) node {\large $1$};
\draw (-0.5,2.8) node {\large $1-\alpha$};
\draw [thick,dotted] (-0,2.8) -- (4.5,2.8);
\draw [black,domain=0:3.82,thin] plot (\x, {2*2.8*exp(-0.008*\x*\x*\x*\x*\x)/(1+exp(-0.00*\x-0.008*\x*\x*\x*\x*\x))});
\draw [black,domain=3.82:4.7,thin] plot (\x, {0.02});
\draw [thin] (0.8,-0.25) -- (0.8,3.55);
\draw (0.9,0.15) node {A};
\draw[fill] (0.8,0) circle (0.03);
\draw [thin] (4,-0.25) -- (4,3.55); 
\draw (4.1,0.15) node {B};
\draw[fill] (4,0) circle (0.03);
\draw (0.4,-0.4) node {{\large $\frac{r_\varepsilon}{\bar r_\varepsilon} \hspace{-0.1cm}\rightarrow \hspace{-0.05cm} 0$}};
\draw (4.6,-0.4) node {{\large $\frac{r_\varepsilon}{\bar r_\varepsilon} \hspace{-0.1cm}\rightarrow \hspace{-0.05cm} +\infty$}};
\draw (2.5,-0.4) node {{\large $\frac{r_\varepsilon}{\bar r_\varepsilon} \hspace{-0.1cm}= \hspace{-0.05cm} \mathcal{O}_\epsilon(1)$}};
\draw [black,domain=0:4.2,thin,dashed] plot (\x, {2*2.8*exp(-0.008*(\x-0.5)*(\x-0.5)*(\x-0.5)*(\x-0.5)*(\x-0.5))/ (1+exp(-0.008*(\x-0.5)*(\x-0.5)*(\x-0.5)*(\x-0.5)*(\x-0.5)))});
\end{tikzpicture}
\end{center}
\caption{\textit{The solid curve represents the function $\Phi(t_{1-\alpha}-u_\epsilon(r_\epsilon))$ while the dashed curve displays the function $\Phi(t_{1-\alpha}-cu_\epsilon(r_\epsilon))$, for some $c\in ]0,1[$ and a radius $r_\varepsilon>0$ satisfying $u_\varepsilon(r_\varepsilon)=\mathcal{O}_\varepsilon(1)$ $($see $($\ref{eq:Fig4-disc}$)$$)$. The solid curve is associated to the sharp asymptotics of Gaussian type for the maximal second kind error probability $\Beta_{\varepsilon,\alpha}(\Theta(r_\varepsilon))$ while the dashed curve is associated to the asymptotic upper bound of the maximal second kind error probability $\Beta_{\varepsilon}(\Theta_a(r_\varepsilon),\Psi_{\tilde D,P})$ of the the spectral cut-off test $\Psi_{\tilde D,P}$, with bandwidth $\tilde D \in \mathbb{N}$ selected as in $($\ref{eq:Dstar_value}$)$.}}
\label{fig:Fig4}
\end{figure}
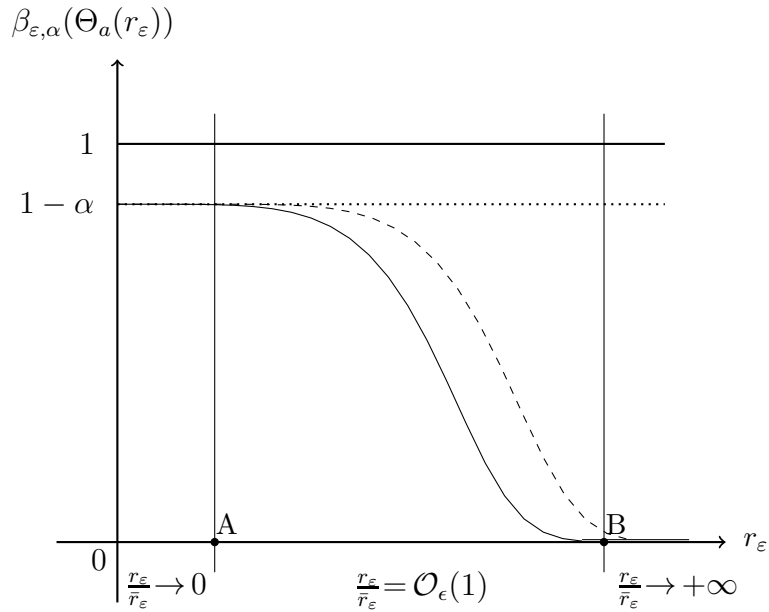

Concerning Proposition \ref{prop:spectral_nonasy}, the following comments are in order:
\begin{itemize}
\item The maximal second kind error probability $\Beta_{\varepsilon}(\Theta_a(r_\varepsilon),\Psi_{\tilde D,P})$ associated to the spectral cut-off test $\Psi_{\tilde D,P}$, with bandwidth $\tilde D \in \mathbb{N}$ selected as in (\ref{eq:Dstar_value}), is asymptotically bounded from above by a quantity that possesses a Gaussian shape. It is worth mentioning that this spectral cut-off test $\Psi_{\tilde D,P}$ is of the same type as the one introduced in Section \ref{subsub:sp-c-off}. In particular, by construction, the spectral cut-off $\Psi_{\tilde D,P}$ is still an $\alpha$-level test. Nevertheless, the bandwidth $\tilde D \in \mathbb{N}$ defined in (\ref{eq:Dstar_value}), is selected in a different manner in order to accommodate the asymptotic paradigm. Indeed, this regularization parameter $\tilde D \in \mathbb{N}$ now depends on the radius $r_\varepsilon$. Notice that this is comparable to the construction of the Ingster test $\Psi_{r_\varepsilon,I}$ introduced in (\ref{yuri_test}), where the Ingster filters $\omega_{j,r_\varepsilon}$ defined in (\ref{yuri_filter}) explicitly depend on the radius $r_\varepsilon$.

\item The asymptotic upper bound of the maximal second kind error probability $\Beta_{\varepsilon}(\Theta_a(r_\varepsilon),\Psi_{\tilde D,P})$ obtained in (\ref{eq:gaussian_shape_cuttoff-0}) is coherent with
the non-asymptotic analysis provided in Section \ref{subsub:sp-c-off} (see, in particular, (\ref{eq:upper_condition2}) and (\ref{bi_nasymp2})). Indeed, in order to guarantee that, for any $\beta \in ]0,1[$,  $\Beta_{\varepsilon}(\Theta_a(r_\varepsilon),\Psi_{\tilde D,P})$ is (asymptotically) upper bounded by $\beta$, we have to solve the equation $a^{-2}_{D} \asymp C_{\alpha,\beta} \varepsilon^2 \sqrt{\sum_{j=1}^{\tilde D} b_j^{-4}}$, for some constant $C_{\alpha,\beta}>0$ (whose value depends on the tools used to control $\Beta_{\varepsilon}(\Theta_a(r_\varepsilon),\Psi_{\tilde D,P})$).

\item In order to conclude our discussion, we provide a heuristic comparison between the asymptotic upper bound of the maximal second kind error probability $\Beta_{\varepsilon}(\Theta_a(r_\varepsilon),\Psi_{\tilde D,P})$ and the sharp asymptotics of Gaussian type obtained in Theorem \ref{thm:1a}. Working as in Section \ref{s:asymp_rate}, we get that, as $\varepsilon \rightarrow 0$,
$$
u_\varepsilon(r_\varepsilon) \sim \frac{r_\varepsilon^4}{\varepsilon^4}\mathcal{R}_\varepsilon^{-(4t+1)} \sim \left(\frac{\mathcal{R}_\varepsilon^{-2s}}{\varepsilon^2 \mathcal{R}_\varepsilon^{2t+1/2}}    \right)^2, \quad \mbox{where} \;\; \mathcal{R}_\varepsilon\;\; \mbox{satisfies}\;\; \mathcal{R}_\varepsilon^{-s} \sim r_\varepsilon.
$$
Note that, thanks to (\ref{eq:comput_0}), $u_\varepsilon(r_\varepsilon) = \mathcal{O}_\varepsilon(1)$ implies that
$r_\varepsilon \sim \varepsilon^{2s/(2s+2t+1/2)}$, as $\varepsilon \rightarrow 0$.
Moreover, in view of (\ref{Fanis-Equiv-Rates-1}), $\tilde D^{-2s} \sim a_{\tilde D}^{-2} \sim r^2_\varepsilon$ as $\varepsilon \rightarrow 0$. Hence, according to the definition of the bandwidth $\tilde D$ given in (\ref{eq:Dstar_value}), as soon as
$u_\varepsilon(r_\varepsilon) = \mathcal{O}_\varepsilon(1)$, in some sense, we have that
$$
u_\varepsilon(r_\varepsilon) \sim \frac{a_{\tilde D}^{-2}}{\varepsilon^2 \sqrt{\sum_{j=1}^{\tilde D} b_j^{-4}}} \quad \mathrm{as} \quad \varepsilon\rightarrow 0.
$$
In particular, it means that we can find a $c \in ]0,1[$ such that
\begin{equation}
\label{eq:gaussian_shape_cuttoff}
\Beta_{\varepsilon}(\Theta_a(r_\varepsilon),\Psi_{\tilde D,P}) \leq \Phi(t_{1-\alpha}-c\,u_\varepsilon(r_\varepsilon)) + o_\varepsilon(1).
\end{equation}
According to Theorem \ref{thm:1a}, it is immediately seen that
\begin{equation}
\label{eq:Fig4-disc}
\Beta_{\varepsilon,\alpha}(\Theta_a(r_\varepsilon))=\Phi(t_{1-\alpha}-\,u_\varepsilon(r_\varepsilon)) + o_\varepsilon(1) < \Phi(t_{1-\alpha}-c\,u_\varepsilon(r_\varepsilon)) + o_\varepsilon(1). 
\end{equation}
Hence, the spectral cut-off test $\Psi_{\tilde D,P}$ defined in (\ref{eq:gaussian_shape_cuttoff-0}), with bandwidth ${\tilde D} \in \mathbb{N}$ selected as in (\ref{eq:Dstar_value}), does not provide sharp asymptotics of Gaussian type. Indeed, it is not designed for that purpose: the spectral cut-off filters associated to this test appear to be quite `rough' in such setting compared to the Ingster filters defined in (\ref{yuri_filter}) (see Figure \ref{fig:Fig4} for a graphical illustration).
%\item According to the discussion above, it seems that is not possible to obtain the optimal constants for the considered minimax signal detection problem with the spectral cut-off test $\Psi_{\tilde D,P}$ defined in (\ref{eq:gaussian_shape_cuttoff-0}), with bandwidth ${\tilde D} \in \mathbb{N}$ selected as in (\ref{eq:Dstar_value}). 
\item If we define a radius $\bar{r}_\varepsilon^\star >0$ to satisfy $u(\bar{r}_\varepsilon^\star)=t_{1-\alpha}-t_\beta$, for prescribed $\alpha, \beta \in ]0,1[$, then, using Theorem \ref{thm:1a}, we immediately get
$$
\Beta_{\varepsilon,\alpha}(\Theta_a(\bar{r}_\varepsilon^\star))=\Phi(t_{1-\alpha}-\,u_\varepsilon(\bar{r}_\varepsilon^\star)) + o_\varepsilon(1)=\beta + o_\varepsilon(1).
$$
Furthermore, according to the definition of the separation radius $\tilde{r}_{\varepsilon, \tilde D}^\star:=r_\varepsilon(\mathcal{E}_a,\Psi_{\tilde D,P},\beta)$ for 
the spectral cut-off test $\Psi_{\tilde D,P}$ defined in (\ref{eq:gaussian_shape_cuttoff-0}), with bandwidth ${\tilde D} \in \mathbb{N}$ selected as in (\ref{eq:Dstar_value}), we have 
$$
\Beta_{\varepsilon,\alpha}(\Theta_a(\tilde{r}_{\varepsilon, \tilde D}^\star) \leq \beta.
$$
However, we conjecture that it is not possible to prove that 
$$
\frac{\tilde{r}_{\varepsilon,\tilde D}^\star}{\bar{r}_\varepsilon^\star}=1+o_\varepsilon(1).
$$
In other words, the spectral cut-off tests appear to be quite `rough' in order to provide the optimal constants of the associated rates for the considered minimax signal detection problem.
\end{itemize}

\begin{remark}
{\rm  Proposition \ref{prop:spectral_nonasy} holds true in a general setting. Indeed, by looking at its proof (the control of $T_1$ and $T_2$), the only condition needed to prove (\ref{eq:gaussian_shape_cuttoff-0}) is that 
\begin{equation}
\label{eq:gen_cond_asy-F}
\exists\; \delta >0 \quad \mbox{such that} \quad \frac{\max_{1\leq j \leq D} b_j^{-2}}{\sqrt{\sum_{j=1}^{D} b_j^{-4}}} = o(D^{-\delta}) \quad \mathrm{as} \quad D \rightarrow +\infty. 
\end{equation}
It is easily seen the condition (\ref{eq:gen_cond_asy-F}) is satisfied in various settings, namely, direct problems (i.e., $b_j =1$, $j\in\mathbb{N}$), well-posed inverse problems 
(i.e., $b_j >b_0$, for some $b_0>0$, $j\in\mathbb{N}$) and mildly ill-posed problems
(i.e., $b_j \asymp j^{-jt}$, $j\in\mathbb{N}$, for some $t>0$).
We point out, however, that it is not satisfied, for instance, in exponential or power-exponential behaviors (i.e., $b_j \asymp e^{-j^rt}$, $j\in\mathbb{N}$, for some $t>0$ and $r\geq1$),
discussed in Remark \ref{rem:2-F}). It is worth mentioning that condition (\ref{eq:gen_cond_asy-F}) is, in general, comparable to the condition $\omega_{0,r_\varepsilon}=o_\varepsilon(1)$, discussed in Theorem \ref{thm:1a}. For more details on the asymptotic expression of $\omega_{0,r_\varepsilon}$ in mildly ill-posed inverse problems, we refer to the proof of Theorem 4.2 of \cite{ISS_2012}, supplementary material, Section 11.3.}
\end{remark}

\section{Conclusions}
\label{s:discussion}
We discussed  non-asymptotic and asymptotic approaches to minimax signal detection trough a unified treatment and provided, in some sense, an overview of this specialized area. % non-asymptotic and asymptotic approaches to minimax signal detection. 
In particular, we considered a Gaussian sequence model that contains classical models as special cases, such as, direct, well-posed inverse and ill-posed inverse problems. We compared the construction of lower and upper bounds for the minimax separation radius (non-asymptotic approach) and the minimax separation rate (asymptotic approach), and brought into light hitherto unknown similarities and links between these two associated minimax signal detection paradigms.  An example of a mildly ill-posed inverse problem was used for illustrative purposes. %In particular, through this example, we demonstrated that non-asymptotic and asymptotic approaches to minimax signal detection, somehow, merge. 
In particular, tools used to derive `asymptotic' results can be exploited to draw `non-asymptotic' conclusions, and vice-versa. To this end, we note that in these considerations we have worked with certain ellipsoids in the space of squared-summable sequences of real numbers, with a ball of positive radius removed, and we confined our attention to the Neyman-Pearson criterion. \\

There are various ways that the above results could be possibly extended. For instance, for the same smoothness classes, similar investigations, could be easily obtained for the total-error probability criterion defined as the sum of the type I and maximal type II error probabilities of a given test $\Psi$, i.e.,
$$
\Zeta_\varepsilon(\Theta_a(r_\varepsilon), \Psi)= \Alpha_\varepsilon(\Psi) + \Beta_\varepsilon(\Theta_a(r_\varepsilon),\Psi),
$$
where $\Alpha_\varepsilon(\Psi)$ and $\Beta_\varepsilon(\Theta_a(r_\varepsilon),\Psi)$ are defined in (\ref{eq:type1-F}) and (\ref{eq:type2}), respectively. Note that, by defining 
$$
\Zeta_\varepsilon(\Theta_a(r_\varepsilon)) = \inf_{\tilde \Psi} [\zeta_\varepsilon(\Theta_a(r_\varepsilon), \tilde \Psi)]
$$
where the infimum is taken over all possible tests $\tilde \Psi$, it is known that (see, e.g., \cite{IS_2003}, Chapter 2) 
that 
$$ \Zeta_\varepsilon(\Theta_a(r_\varepsilon)) = \inf_{\alpha \in ]0,1[} \left[ \alpha + \Beta_ {\varepsilon,\alpha}(\Theta_a(r_\varepsilon)) \right],
$$
where $\Beta_ {\varepsilon,\alpha}(\Theta_a(r_\varepsilon))$ is the minimax second kind error probability defined in Definition \ref{mskep-F}.\\

Similar investigations for the  Neyman-Pearson criterion and/or the total-error probability criterion should also be possible for other classes ${\cal F}$ of signals, such as those characterized by their non-zero coefficients (dense or sparse signals) and $l_p$-bodies with $p \in ]0, 2]$ (see, e.g.,  \cite{S_1996}, \cite{Baraud}, \cite{ISS_2012}, \cite{LLM_2012}). In the same spirit, several contributions have been proposed in various regression and density models which provide attractive frameworks for investigation in the minimax testing theory (see, e.g., \cite{GL_2002}, \cite{E_1994}, \cite{FL_2006}, \cite{Butucea}, \cite{IS_2009}, \cite{BMP_2009}, \cite{LPN_2014}). \\

For the sake of brevity and clarity in our presentation, we have also not discussed adaptation issues of the involved testing procedures in the considered minimax signal detection paradigms. Indeed, the filters used to design the spectral cut-off (non-asymptotic framework) and Ingster (asymptotic framework) tests explicitly depend on the form of the sequence  $(a_j)_{j\in\mathbb{N}}$ that measures the smoothness of the signal $\theta$, which is, in general, unknown in practice. It is therefore of paramount importance in practical applications to provide minimax testing procedures that do not explicitly depend on the associated smoothness parameter.  This is, usually, referred to as the `adaptation' problem (see, e.g., \cite{Baraud}, \cite{IS_2003}, \cite{ISS_2012}, \cite{MM_2013}).\\

However, all the above investigations need careful attention that is beyond the scope of the present work.

\bigskip

\bibliography{Survey}
\bibliographystyle{plain}

\end{document}